\documentclass[sn-mathphys]{sn-jnl}% Math and Physical Sciences Reference Style
\usepackage{amsthm}
\usepackage{comment}
\usepackage{mathrsfs}

\jyear{2025}%

\theoremstyle{thmstyleone}%
\newtheorem{theorem}{Theorem}[section]% meant for sectionwise numbers
\newtheorem{proposition}[theorem]{Proposition}%
\newtheorem{lemma}[theorem]{Lemma}%
\newtheorem{corollary}[theorem]{Corollary}%
\newtheorem{remark}[theorem]{Remark}%
\newtheorem{definition}[theorem]{Definition}%

\theoremstyle{thmstyletwo}%

\renewcommand{\proofname}{Proof.}

\usepackage{anyfontsize}
\usepackage{enumitem}

\newcommand{\Integer}{ \mathbb{Z} }

\newcommand{\Real}{ \mathbb{R} }

\newcommand{\ie}{{\it i.e.}}

\newcommand{\tnorm}[1]{\lvert\!\lvert\!\lvert #1 \rvert\!\rvert\!\rvert}

\newcommand{\red}[1]{\textcolor{black}{#1}}

\begin{document}
%\pagewiselinenumbers %To be deleted

\title[Well-posedness of combustion problem]{Well-Posedness for the Cauchy Problem with Free Boundary in Half-Space Arising from Combustion Problem}


\author*[1]{\fnm{} \fnm{Ken} \sur{Furukawa}} \email{furukawa@sci.u-toyama.ac.jp}
\author[2]{\fnm{} \fnm{Yoshikazu} \sur{Giga}} \email{labgiga@ms.u-tokyo.ac.jp}
\author[3]{\fnm{} \fnm{Naoto} \sur{Kajiwara}} \email{kajiwara.naoto.p4@f.gifu-u.ac.jp}



\affil*[1]{\orgdiv{Faculty of Science, Academic Assembly}, \orgname{University of Toyama}, \orgaddress{\street{Gofuku 3190}, \city{Toyama}, \postcode{930-0887}, \state{Toyama}, \country{Japan}}}
\affil[2]{\orgdiv{Graduate School of Mathematical Sciences}, \orgname{The University of Tokyo} \orgaddress{\street{3-8-1 Komaba, Meguro-ku}, \postcode{153-8914}, \state{Tokyo}, \country{Japan}}}
\affil[3]{\orgdiv{Applied Physics Course, Department of Electrical, Electronic and Computer Engineering}, \orgname{Gifu University}, \orgaddress{\street{1-1 Yanagido}, \city{Gifu}, \postcode{501-1193}, \state{Gifu}, \country{Japan}}}
%\affil[3]{\orgdiv{Department}, \orgname{Organization}, \orgaddress{\street{Street}, \city{City}, \postcode{610101}, \state{State}, \country{Country}}}


\abstract{
In this paper, we address the local-in-time well-posedness of an obstacle problem arising from a combustion model.
The equation under consideration is a linear inhomogeneous heat equation on the moving domain $(s(t),\infty)$, subject to homogeneous Dirichlet and inhomogeneous Neumann boundary conditions.
Our main result establishes the existence and uniqueness of a positive classical solution.

Although several results are known for the homogeneous heat equation, much less is known about the inhomogeneous case.
We show that the boundary function $s(t)$ is $\frac{1}{2}$-H\"older continuous in time when the forcing term has a singularity of order $O(1/\sqrt{t})$ as $t \to 0$.
We also prove that if both the initial data and the forcing term are smooth, then $s$ has higher regularity than $\frac{1}{2}$-H\"older continuity.
}

%Data assimilation, the primitive equations, maximal regularity



\maketitle

%\backmatter
%\bmhead{Acknowledgments}

%Remark for obstacle problem

\section{Introduction} \label{intro}
\begin{figure}[b]
    \centering
    \includegraphics[width=0.6\textwidth]{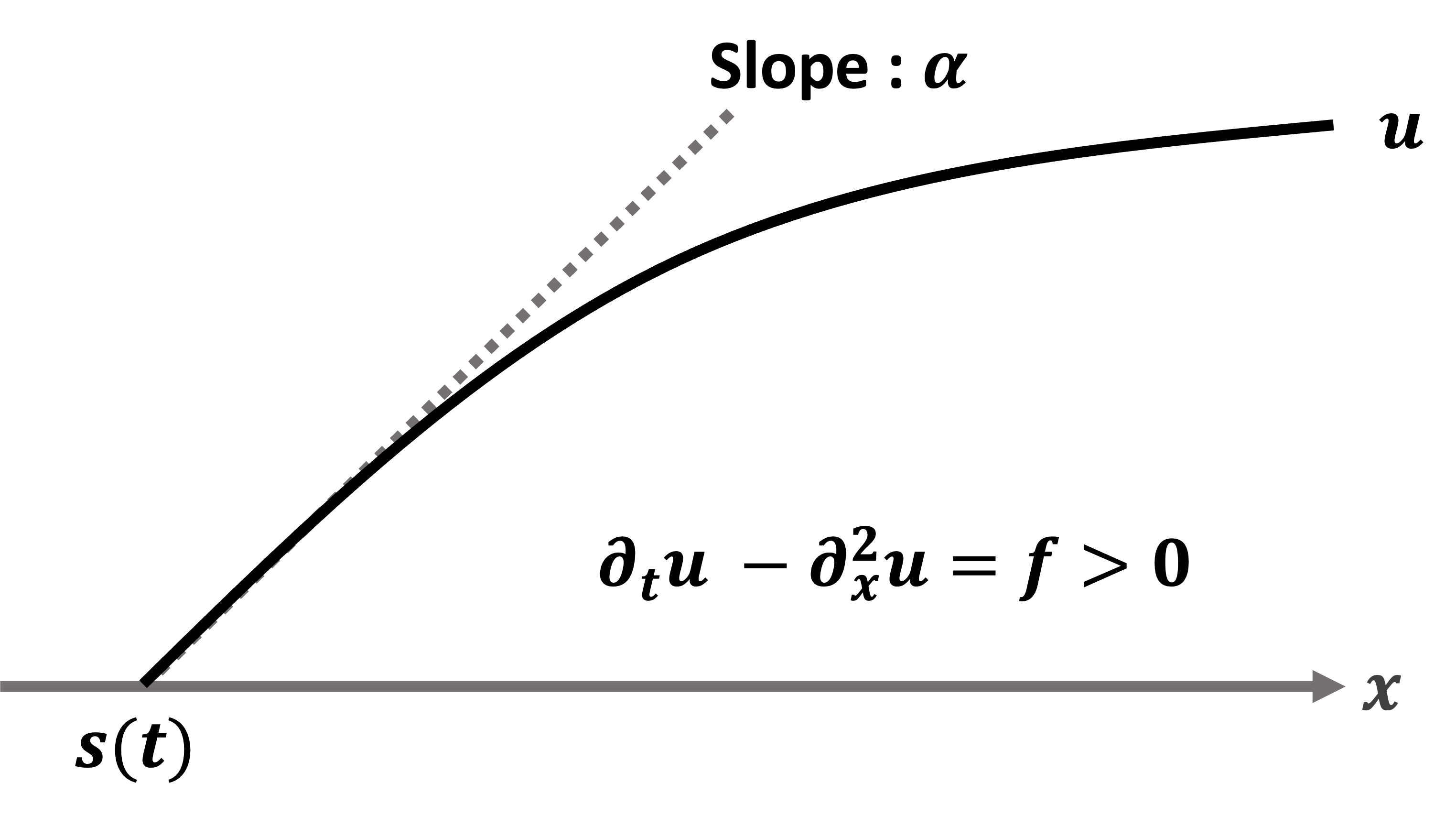}
    \caption{Profile of $u(\cdot,t)$ with positive forcing. The free boundary is located at $x=s(t)$, where $u=0$ and the boundary slope is $\alpha$.}
    \label{fig:u-profile}
\end{figure}

We consider the well-posedness of the heat equation with the free-boundary problem
\begin{equation} \label{eq_original_fb_problem}
    \begin{aligned}
        \partial_t u - \partial_x^2 u
        & = f
        & \quad \text{in}
        & \quad Q_T, \\
        - \frac{\partial u}{\partial \nu}
        & = \alpha
        & \quad \text{on}
        & \quad \Gamma_T, \\
        u
        & = 0
        & \quad \text{on}
        & \quad \Gamma_T, \\
        u
        & \geq 0
        & \quad \text{in}
        & \quad Q_T, \\
        u \vert_{t=0}
        & = u_0
        & \quad \text{in}
        & \quad \Omega_0.
    \end{aligned}
\end{equation}
Here $\nu$ is the outward normal, $\Omega_t = \{ x \in \Real \,;\, u(x, t) > 0 \}$, $Q_T = \{ (x, t) \in \Real \times (0, T) \,;\, u(x, t) > 0 \}$, and $\Gamma_T = \{ (x, t) \in \Real \times (0, T) \,;\, u(x, t) = 0 \}$.
We denote by $\Omega = (0, + \infty)$ the initial domain.
Figure~\ref*{fig:u-profile} illustrates a schematic profile of $u(\cdot,t)$ at a fixed time $t$ in the one-dimensional setting $\Omega_t=(s(t),\infty)$.
%At the free boundary $x=s(t)$, the solution vanishes, while the prescribed flux condition gives $\partial_x u(s(t),t)=\alpha$, since the outward normal is $\nu=-1$.
%The dotted line represents the tangent to the profile at the free boundary, with slope $\alpha$; in the positive region $x>s(t)$, the solution satisfies $\partial_t u-\partial_x^2 u=f$, with $f>0$ in the case illustrated.}
%\begin{figure}[htbp]
When the initial data $u_0$ and the forcing term $f$ are positive, the system can describe a combustion problem.
The solution $u$ describes the flame propagation front.
In combustion-type free boundary problems, the boundary condition of constant normal heat flux originates from the thin flame approximation.
For readers interested in the combustion-theoretic background of this model, we refer to the book by Williams \cite{Williams1985}.
The reader is also referred to \cite{GalaktionovHulshofVazquez1997} for \red{the mathematical and historical background}.
In contrast to Stefan-type problems, the second equation of \eqref{eq_original_fb_problem} does not describe the propagation speed of the interface directly.
%The boundary conditions are overdetermined; nevertheless, since the interface itself moves, the system is expected to be well-posed within an appropriately chosen mathematical framework.

In the homogeneous case, namely when $u$ solves the heat equation on $Q_T$, several results are available concerning existence and asymptotic behavior.
Caffarelli and V\'{a}zquez \cite{CaffarelliVazquez1995} proved the existence of a weak solution with a compactly propagating front that vanishes in finite time by means of variational techniques.
Hilhorst and Hulshof \cite{HilhorstHulshof1991} studied a one-dimensional travelling front problem with a positive heat source at the front.
In \cite{HilhorstHulshof1994}, they also showed the asymptotic profile of the solution at the vanishing time, imposing a homogeneous Neumann boundary condition at the left endpoint of the interval.
Caffarelli, Lederman, and Wolanski \cite{CaffarelliLedermanWolanski1997uniform,CaffarelliLedermanWolanski1997pointwise} studied a two-phase problem, establishing the existence of weak solutions and describing their asymptotic behavior.
Weiss \cite{Weiss2003,Weiss2004} investigated the profile of the free boundary $\{u = 0\}$.
Choi \cite{Choi2009} studied a related problem on a cylindrical domain and established results on the behavior and regularity of the solution at contact points.
For a more detailed overview of the history and previous studies on combustion problems, we refer the reader to \cite{FurukawaGigaKajiwara2026}.
\red{The reader is also referred to \cite{Friedman1964, ProtterWeinberger1967, LadyzenskajaSolonnikovUralceva1968,Hulshof1987, Lunardi1995, Vazquez1996, Kim2003, PucciSerrin2007, PetrosyanYip2008} for related studies on free boundary problems and heat equations.}
%LadyzenskajaSolonnikovUralceva1968Hulshof1987
For inhomogeneous cases, the authors \cite{FurukawaGigaKajiwara2025} proved the nonexistence of bounded solutions with zero initial data and bounded positive external forcing in the spatial domains $(0, \infty)$ and $d$-dimensional bent half-spaces.
They also established the unique existence of \red{a} bounded self-similar solution with zero initial data and spatially constant forcing of order $\frac{1}{\sqrt{t}}$.
%Such solutions are necessarily unbounded.

In this paper, the interface is assumed to be located at $x = s(t)$.
The initial interface may be located at the origin, \ie, $s(0) = 0$.
Under this assumption, the third equation can be written as
\begin{align*}
    u(s(t), t)
    \equiv 0,
    \quad t \geq 0.
\end{align*}
This identity leads to the ODE governing the motion of $s$:
\begin{align*}
    \alpha \dot{s} (t)
    + \partial_x^2 u (s(t), t)
    + f(s(t), t)
    = 0,
    \quad t>0.
\end{align*}
Here we used the equation $\partial_t u = \partial_x^2 u + f$.
We introduce the shifted unknown
\begin{align*}% \label{eq_definition_of_v_g_g0}
    \begin{split}
        v(y, t)
        := u(y + s(t), t)
    \end{split}
\end{align*}
for $y \in \Real$.
The pair of functions $(v, s)$ then satisfies
\begin{equation} \label{eq_v}
    \begin{aligned}
        \partial_t v - \partial_y^2 v
        & =  \dot{s} \partial_y v
        + f(\cdot + s, t),
        & y \in \Omega, \, t > 0, \\
        \alpha \dot{s} 
        & = - \gamma_0 \partial_y^2 v
        - f(s, t),
        & t >0, \\
        \gamma_0 \partial_y v
        & = \alpha,
        & t > 0,\\
        v
        & = v_0, \,
        s
        = 0
        & t = 0,
    \end{aligned}
\end{equation}
where $\gamma_0$ denotes the trace operator at $y = 0$.

In order to state the notion of the solution, we introduce the function space $X^\theta(T)$, which is a subspace of the $\theta$-H\"{o}lder space, defined by
\begin{align*}% \label{eq_condition_for_boundary_function}
    \begin{split}
        X^\theta(T)
        & = \biggl \{
            s
            \in C[0, T]
            \cap C^1(0, T]
            \; ; \;
            s(0)
            = 0, \;\\
        & \quad\quad\quad\quad\quad\quad
            \Vert
                s
            \Vert_{X^\theta(T)}
            = \sup_{0 < \tau < T}
                \tau^{1 - \theta}
                \vert
                    \dot{s}(\tau)
                \vert
            < \infty
        \biggr \}.
    \end{split}
\end{align*}
We also define $Y(T)$ as the space of $C^1$ functions on $[0, T]$ satisfying $s(0) = 0$, equipped with the norm $\Vert s \Vert_{Y(T)} = \Vert \dot{s} \Vert_\infty$, where $\Vert \cdot \Vert_{\infty}$ denotes the uniform norm.
Clearly, $Y(T) \subset X^\theta(T)$.
For $q \in [1, \infty]$ and $m \in \Integer_{\geq 0}$, we denote by $L^q(D)$ the Lebesgue space and by $W^{m, q}(D)$ the Sobolev space on a domain $D$.
We denote by $BUC(D)$ the space of bounded and uniformly continuous functions on $D$.
We denote by $BUC^m(D)$ the space of $m$-times differentiable functions whose derivatives up to order $m$ belong to $BUC(D)$ for $m \in \Integer_{\geq 0}$.
%Analogously, we denote by $BC(D)$ and $BC^m(D)$ the spaces of bounded continuous functions and bounded $m$-times continuously differentiable functions, respectively.
The natural scaling for solutions to the heat equation satisfying $\gamma_0 v = 0$ and $\gamma_0 \partial_y v = \alpha$ is
\begin{align*}
    \sqrt{\lambda}
    v
    \left(
        \frac{x}{\sqrt{\lambda}},
        \frac{t}{\lambda}
    \right), \;
    \sqrt{\lambda} s
    \left(
        \frac{t}{\lambda}
    \right) \quad
    \lambda > 0.
\end{align*}
The corresponding scale-critical space for $s$ is $X^\frac{1}{2}(T)$.

Assume that the external force $f$ satisfies the following conditions:
\begin{enumerate}[label=(\Alph*).]
    \item The forcing $f \in C((0, T); BUC^2(\Real))$ is a positive function satisfying
    \begin{gather*}
        t^\frac{1}{2} f, \,
        t^{1 - \varepsilon} \partial_x f,\,
        t^{\frac{3}{2} - \varepsilon} \partial_x^2 f
        \in BUC(\Real \times (0, T)),\\
        t^{1 - \varepsilon} \partial_x f,\,
        t^{\frac{3}{2} - \varepsilon} \partial_x^2 f
        \in BUC((0, T); L^1(\Real)).
    \end{gather*}
    for some $\varepsilon \in (0, \frac{1}{2})$.
    \item The forcing term $f$ is parabolically H\"older continuous on $\Real \times (0,T)$ in the sense that
    \begin{align*}
        C_f
        & := \sup_{0 < \tau < t < T}
            \sqrt{t \tau}
                \sup_{x_1, x_2 \in \Real, \, x_1 \neq x_2}
                    \frac{
                        \vert
                            f(x_2, t)
                            - f(x_1, \tau)
                        \vert
                    }{
                        \vert
                            x_2
                            - x_1
                        \vert
                        + \vert
                            t
                            - \tau
                        \vert^\frac{1}{2}
                    }
        & < \infty.
    \end{align*}
    \item There \red{exist} $r_0 > 0$ and $t_0 > 0$ such that
    \begin{align*}
        \inf_{t \in (0, t_0)}
            \inf_{x \in (-r_0, r_0)}
                t^\frac{1}{2}
                f(x, t)
        >0.
    \end{align*}
    
\end{enumerate}
These three conditions capture singular \red{solutions}.
The first \red{result} of this paper is
\begin{theorem} \label{thm_main1}
    Let $T > 0$ and $\alpha > 0$.
    Assume (A)-(C).
    Let $v_0 \in BUC^1(\Omega)$ be a positive function satisfying $\partial_y v_0 \in L^1(\Omega)$, $\gamma_0 v_0 = 0$, and $\gamma_0 \partial_y v_0 = \alpha$.
    Then there exist $\varepsilon > 0$ and $T_\ast \in (0, T)$ \red{such that}, if
    \begin{align*}
        C_f
        + \sup_{0 < t < T}
            t^\frac{1}{2}
            \Vert
                f(t)
            \Vert_{BUC(\Real)}
        < \varepsilon,
    \end{align*}
    there exists a unique solution
    \begin{gather*}
        v
        \in C([0, T_\ast]; BUC(\Omega))
        \cap C^1((0, T_\ast]; C(\Omega))
        \cap C((0, T_\ast]; BUC^1(\Omega) \cap C^2(\Omega)), \\
        v(t) > 0 \text{ on $\Omega$} ,\\
        s
        \in X^\frac{1}{2}(T_\ast),
    \end{gather*}
    to \eqref{eq_v} associated with the boundary condition
    \begin{align} \label{eq_bc_in_theorem}
        \gamma_0 v
        = 0.
    \end{align}
\end{theorem}

\begin{remark}
    A spatially constant forcing of the form
    \begin{align*}
        f(x, t)
        = \frac{f_0}{\sqrt{t}}
    \end{align*}
    for some small constant $f_0 > 0$ satisfies the assumptions of Theorem \ref{thm_main1}.
\end{remark}
\begin{remark}
    Theorem \ref{thm_main1} means that to obtain the solution in \red{the scale-critical} space $s \in X^\frac{1}{2}(T)$, \red{a} singularity of $f$ is needed.
    If $f$ is smoother than $O(\frac{1}{\sqrt{t}})$, a \red{scale-critical} solution \red{does not seem to} exist; see also Remark \ref{rmk_vanishing_boundary_layer}.
\end{remark}

\begin{remark}
    Theorem \ref{thm_main1} is not inconsistent with our previous result (Theorem 7 in \cite{FurukawaGigaKajiwara2025}).
    In our previous work, we proved that, when the initial datum is identically zero and the external forcing is given by $f=h/\sqrt{t}$ with $h>0$, the interface $s(t)$ moves to the left if $h>\alpha/\sqrt{\pi}$, whereas it moves to the right if $h<\alpha/\sqrt{\pi}$.
    The crucial difference between our previous work and the present study is that no initial interface exists in the former, whereas an initial interface is already present in the latter. 
    Consequently, when an initial interface exists and the external forcing (heat input) has a singularity of order $1/\sqrt{t}$, the contact point moves to the left; that is, the interface expands.
    The smallness assumption on $f$ in Theorem \ref{thm_main1} is imposed solely to apply the fixed-point theorem, and it is of a different nature from the condition $h<\alpha/\sqrt{\pi}$ in our previous work.
    The latter condition reflects the competition between the singular external forcing and the heat absorption at the contact point.
    Although the singular forcing immediately induces the formation of an interface, the heat absorption at the contact point dominates when $h<\alpha/\sqrt{\pi}$, and consequently the interface $s(t)$ moves to the right, namely in the direction of shrinking.
    In contrast, since an initial interface has already been formed in the present setting, it is natural to expect that, if $f$ is sufficiently singular, the interface will instead move in the direction of expansion.
\end{remark}

The second result is the existence of \red{a} smooth solution.
In the following theorem, we do not assume (A)-(C).
\begin{theorem} \label{thm_main2}
    Let $T > 0$ and $\alpha > 0$.
    Let $v_0 \in BUC^2(\Omega)$ be a positive function satisfying $\partial_y v_0 \in L^1(\Omega)$, $\gamma_0 v_0 = 0$, and $\gamma_0 \partial_y v_0 = \alpha$.
    Let $f \in BUC((0, T); BUC^2(\Real)) \cap Lip(\Real \times (0, T))$ be a positive function.
    Assume further
    \begin{align} \label{eq_compatibility_thm2}
        \partial_y^2 v_0 (0)
        + f(0,0)
        \neq 0.
    \end{align}
    Then there exist $T_\ast \in (0, T)$ and a unique solution
    \begin{gather*}
        v
        \in C([0, T_\ast]; BUC^2(\Omega))
        \cap C^1((0, T_\ast]; BUC(\Omega)),\\
        v(t) > 0 \text{ on $\Omega$} ,\\
        s
        \in Y(T_\ast)
    \end{gather*}
    to \eqref{eq_v} associated with the boundary condition \eqref{eq_bc_in_theorem}.
\end{theorem}
\begin{remark}
    The condition \eqref{eq_compatibility_thm2} ensures that the contact point $s(t)$ starts moving immediately from its initial position.
    Similarly, an analogous assumption is imposed in the form of {\rm(H1)-(H2)} in \cite{CaffarelliVazquez1995}.
    
    Moreover, assumption {\rm (H1)} imposes a monotonicity condition on the initial data, whereas Theorem \ref{thm_main2} does not require any such geometric assumption on the initial data.

    In \cite{HilhorstHulshof1991}, the initial data are assumed to decay at infinity, and the existence of a global-in-time solution is established in the presence of a positive heat flux $q>0$ at the origin and in the absence of external forcing.
    In contrast, in Theorem \ref{thm_main2}, we consider the case $q=0$ without assuming any decay of the initial data at infinity.
    However, global-in-time extendability of the local solution is not addressed.
\end{remark}

We introduce the strategy for proving our results.
To control the highest order term $- \gamma_0 \partial_y^2 v$ in \eqref{eq_v}, one would need regularity higher than the natural regularity of solutions provided by the first equation.
However, it seems difficult to find suitable Sobolev spaces as solution spaces for the drift-diffusion equation that are compatible with the two boundary conditions.
Nevertheless, there is a more reasonable approach.
We therefore consider the derivative of $v$.
We write
\begin{gather} \label{eq_definition_w_w0_h}
    w
    = \partial_y v, \quad
    w_0
    = \partial_y v_0.
\end{gather}
Then $(w, s)$ is the solution to
\begin{equation} \label{eq_w}
    \begin{aligned}
        \partial_t w - \partial_y^2 w - \dot{s} \partial_y w
        & = \partial_y f(\cdot + s, \cdot),
        & y \in \Omega, \, t>0, \\
        \alpha \dot{s}
        & = - \gamma_0 \partial_y w
        - f(s, \cdot),
        & t>0, \\
        w
        & = \alpha>0,
        & y=0, \, t>0, \\
        w
        & = w_0,
        & y \in \Omega, \, t=0, \\
        s
        & = 0,
        & t=0.
    \end{aligned}
\end{equation}
The solution $v$ is recovered by the fundamental theorem of calculus
\begin{align} \label{eq_v_frome_w}
    v(y, t)
    := \int_{0}^y
        w(z, t)
    dz,
\end{align}
if this integral makes sense for all $t > 0$ and $y \in \Omega$.
\begin{remark} \label{rmk_boundary_condition}
    This formulation seems natural, since the boundary conditions for $v$ are automatically satisfied.
    It is worth mentioning that $v$ satisfies the first equation in \eqref{eq_v} only when $s$ satisfies the second equation of \eqref{eq_w}.
    Indeed,
    \begin{align*}
        & \partial_t v(y,t)\\
        %& = \int_{0}^y
        %    \partial_t w (y, t)
        %dy\\
        & = \int_{0}^y
            \partial_y^2 w
            + \dot{s}(t) \partial_y w
            + \partial_y f(y + s(t), t)
        dy\\
        & = (
            \partial_y w(y, t)
            - \partial_y w(0, t)
        )
        + \dot{s}(t) (
            w(y, t)
            - \alpha
        )
        + (
            f(y + s(t), t)
            - f(s, t)
        )\\
    %    & = \partial_y w(y, t)
    %    + \dot{s}(t) w(y, t)
    %    + f(y + s(t), t).
        & = \partial_y^2 v(y, t)
        + \dot{s}(t) \partial_y v(y, t)
        + f(y + s(t), t).
    \end{align*}
    If $s$ does not satisfy \eqref{eq_v}, then singular terms appear from the boundary values in the equation for $v$.
\end{remark}
The boundary term $\gamma_0 \partial_y w$ in \eqref{eq_w} is relatively easier to control by the natural regularity from the first equation of \eqref{eq_w}.
The system of equations \eqref{eq_w} involves only a single boundary condition, and therefore it is considerably easier to analyze than \eqref{eq_v}.
To prove the theorem, we consider the mild solution $(w, s)$ \red{defined as follows.}
\begin{definition}
    Let $T>0$.
    Let
    \begin{align} \label{eq_assumptions_w0_f_in_def_of_mild_solution}
        w_0 \in BUC(\Omega), \quad
        %f \in C((0, T); BUC^2(\Omega))
        f \in C((0, T); BUC^{1}(\Omega))
    \end{align}
    We call the pair of functions
    \begin{align*}
        (w, s)
        \in C([0, T); BUC(\Omega))
        \cap C((0, T); BUC^1(\Omega))
        \times C[0, T)
    \end{align*}
    a mild solution of \eqref{eq_w} when, for all $t \in (0, T)$,
    \begin{align} \label{eq_mild_solution_w}
        \begin{split}
            & w(y, t)\\
            & = \int_0^\infty
                G_t (y - z + s(t)) w_0(z)
            dz\\
            & + \int_0^t
                G_{t - \tau} (y + s(t) - s(\tau)) f(s(\tau), \tau)
            d\tau\\
            & - \int_0^t
                \alpha \partial_y G_{t - \tau}(y + s(t) - s(\tau))
            d\tau\\
            & + \int_0^t
                \int_{0}^\infty
                    G_{t - \tau}(y - z + s(t) - s(\tau)) \partial_z f(z + s(\tau), \tau)
                dz
            d\tau
        \end{split}
    \end{align}
    and
    \begin{align} \label{eq_mild_solution_s}
        s(t)
        = - \int_0^t
            \frac{1}{\alpha}
            \left(
            \gamma_0 \partial_y w (\tau)
                + f(s(\tau), \tau)
            \right)
        d\tau
    \end{align}
    hold, where $G_t(x) = \frac{1}{\sqrt{4 \pi t}} e^{-\frac{x^2}{4t}}$ denotes the one-dimensional Gaussian kernel.
\end{definition}
It is worth \red{mentioning} that the mild solution formula and the Young inequality provide \red{an} $L^1(\Omega)$-estimate for $w$ under suitable integrability conditions for $w_0$ and $f$, which leads to \red{the} boundedness of $v$.

%\begin{figure}[htbp]
\begin{figure}[t]
    \centering
    \includegraphics[width=\textwidth]{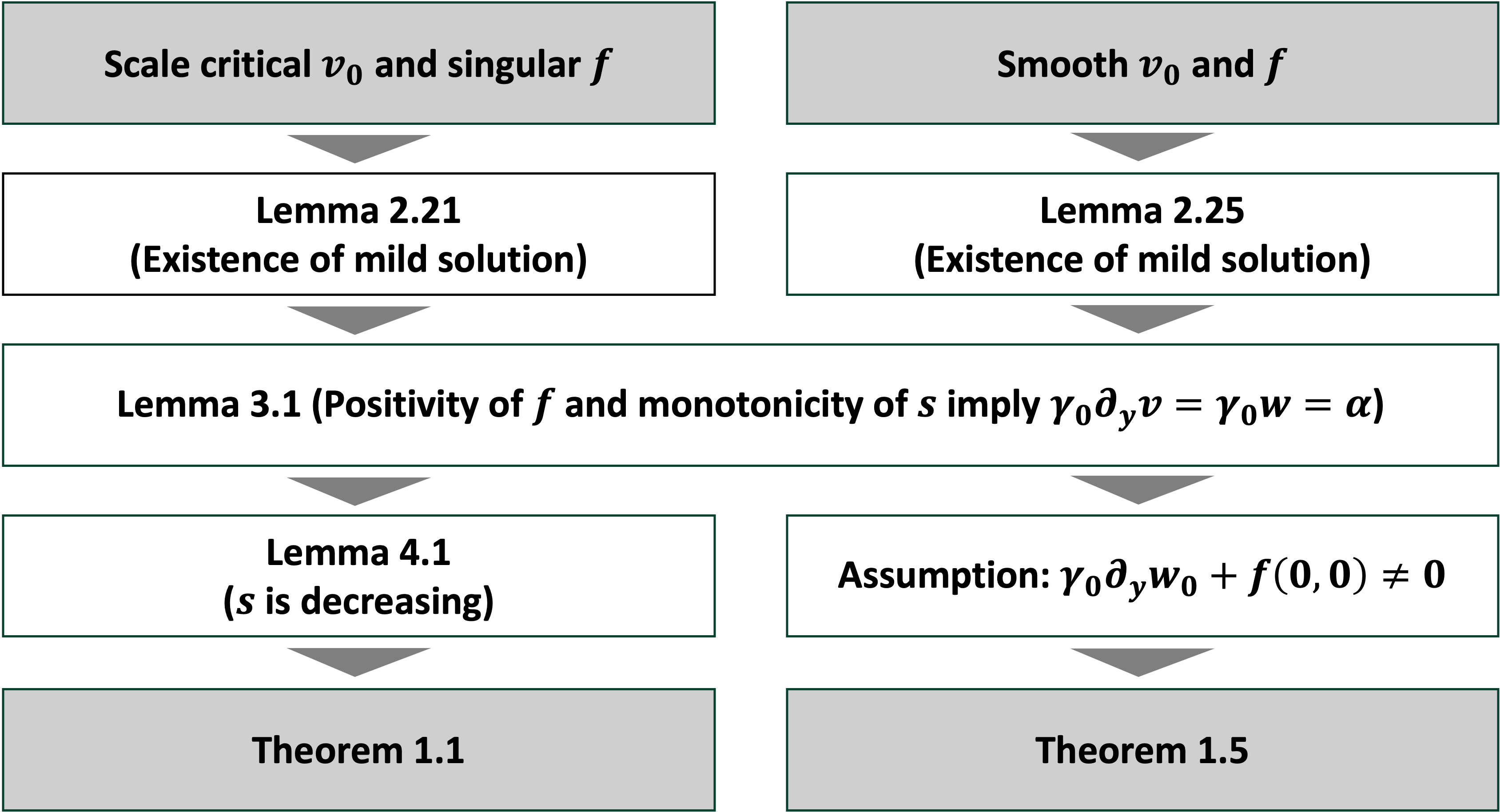}
    \caption{Outline of the proof strategy for singular forcing (left) and smooth data (right).}
    \label{fig:proof-outline}
\end{figure}

We derive this formula \eqref{eq_mild_solution_w} from the equation by the Laplace transformation\red{,} but we \red{note that} this formula \red{can} also \red{be} obtained by Green's formula based on integration by parts.
The proofs of Theorems~\ref*{thm_main1} and~\ref*{thm_main2} follow the three-step scheme in Figure~\ref*{fig:proof-outline}.
They rely on Young's inequality for convolution estimates and do not require general maximal regularity theory as \red{in} \cite{FurukawaGigaKajiwara2026}.
The first step is the construction of the mild solution $(w, s)$ without yet verifying the boundary condition $\gamma_0 w = \alpha$.
The appropriate function space is associated with the mixed norm
\begin{align} \label{eq_norm_for_w_in_introduction}
    \sup_{0 < t < T}
        \Vert
            w(t)
        \Vert_{\infty}
    + \sup_{0 < t < T}
        t^\frac{1}{2}
        \Vert
            \partial_y w(t)
        \Vert_{\infty},
\end{align}
which is quite natural from the viewpoint of parabolic scaling.
The basic idea is to estimate the integral kernels in \eqref{eq_mild_solution_w} and derive a quadratic bound for $\partial_y w$ in terms of $\Vert s \Vert_{X^\frac{1}{2}(T)}$.
The solution $w$ of \eqref{eq_mild_solution_w} is determined from $s$ and does not contain any nonlinearity by itself.
Thus the main task in this step is the construction of $s$ in \eqref{eq_mild_solution_s}, regarding $\gamma_0 \partial_y w$ as a nonlinear function of $s$.
Since the scale of self-similarity is given by $s(t) \propto \sqrt{t}$ (see also \cite{FurukawaGigaKajiwara2025} for self-similar solutions), the base space for $w$ must be a class of bounded continuous functions in order for the term $\gamma_0 \partial_y w$ to be well defined.
To ensure the boundedness of $v$, we also require the integrability of $w$.
This integrability is proved in essentially the same way as in the proof of estimating \eqref{eq_norm_for_w_in_introduction}, replacing the uniform norm with the $L^1$-norm and assuming suitable integrability of the initial data $w_0$ and the forcing term $f$.
To estimate $\gamma_0 \partial_y w$, the first and third terms on the right-hand side of \eqref{eq_mild_solution_w}, which involve the boundary layer effect, need to be treated carefully.
We derive the formulas \eqref{eq_formula_w_1_w_3}--\eqref{eq_formula_diff_w_1_w_3}, in which the boundary-layer effects are moderated, making the estimates more tractable, and then estimate these terms in a unified manner, since they cannot be controlled separately.

The boundary condition is verified in Lemma \ref{lem_Neumann_boundary_condition} under the assumption that $s$ is monotone.
This assumption is natural because, in general, overdetermined problems on fixed domains admit no solution.
In this section, we derive the integral equation satisfied by $w$ from the ODE satisfied by $s$, and compare it with \eqref{eq_mild_solution_w}.
This is done by multiplying \red{the second equation of \eqref{eq_w} by} $G_{t - \tau}(y + s(t) - s(\tau))$ and integrating over the parabolic cylinder $\Omega \times (0, T)$.
Using integration by parts and comparing the resulting integral formulation with that of $w$, we find that an integral involving $\gamma_0 w - \alpha$ multiplied by a function with a consistent sign must vanish.
This leads to the boundary condition for $w$.

In order to determine the sign of $\dot{s}$ for $s \in X^\frac{1}{2}(T)$ on a short interval around $t=0$, we consider the sign of $\sqrt{t}(\gamma_0 \partial_y w + f(s, \cdot))$.
The formula \eqref{eq_gamma_0_partial_y_w} in Section \ref{sec_moving_boundary_condition} is crucial for understanding the behavior of $\gamma_0 \partial_y w$.
Since the integral $\sqrt{t} \int_0^\infty \partial_y G_t (- z + s(t)) (w_0(z) - \alpha) dz$ and the contributions of the forcing term are $o(1)$, using the formula \eqref{eq_gamma_0_partial_y_w} we observe that $\sqrt{t} \gamma_0 \partial_y w$ can be small near $t = 0$.
Therefore, if $f \propto \frac{1}{\sqrt{t}}$ near $y = 0$, as in Assumption (C), we obtain $\dot{s} < 0$.

The proof of Theorem \ref{thm_main2} is essentially based on the same scenario as the proof of Theorem \ref{thm_main1}.
The main differences lie in the solution space to which the mild solution belongs and in the method used to determine the sign of $\dot{s}$.
Regarding the latter, since we assume \eqref{eq_compatibility_thm2}, the continuity of $\gamma_0 \partial_y w$ and $f(s,t)$ implies that the sign of $\dot{s}$ remains consistent on some short time interval.
Since $\partial_y w_0$ is uniformly continuous, the singularity from the boundary layer effect \red{fortunately vanishes}; see the formula \eqref{eq_formula_w_1_w_3_when_w0_is_differentiable} and Remark \ref{rmk_vanishing_boundary_layer}.
Consequently, we can construct $s$ in $Y(T)$ using Banach's fixed point theorem, and then obtain $w \in BUC^1(\Omega)$ on $[0, T_\ast]$.
Even if we assume that $w_0 \in BUC^1(\Omega)$, the integral involving the initial data in the mild formulation for $\partial_y w$ produces a Gaussian term arising from the boundary values.
Since this Gaussian term is singular of order $O(t^{-1/2})$, it may at first sight seem impossible to obtain a uniform estimate for $\partial_y w$ in $t$.
However, the derivative of the third term on the right-hand side of the mild formulation of $w$ cancels this Gaussian contribution.
Note that the boundary condition for $w_0$ is essentially used here.

We use the following notation.
%We denote by $\Vert f \Vert_{\infty}$ the uniform norm of a continuous function $f$.
%We denote by $BUC^m(\Omega) \subset C^m(\Omega)$ the space of $m$-times continuously differentiable functions whose $l$-th derivatives for $l \leq m$ are uniformly continuous and bounded.
%xFor a smooth function $\varphi$, we write $\varphi_t(x)$ to denote the parabolically scaled function $\frac{1}{\sqrt{t}}\varphi(\frac{x}{\sqrt{t}})$.
We denote by $r^\prime = \frac{r}{r - 1}$ the H\"{o}lder conjugate index for $r \in [1, \infty]$, with the standard modification for $r=1, \infty$.

This paper is organized as follows.
In Section \ref{sec_mild_solution}, we first construct the mild solution \eqref{eq_mild_solution_w}.
In Section \ref{sec_boundary_condition}, we prove that when $s$ is monotone, \ie, when the front is moving, the solution $w$ satisfies the boundary condition $\gamma_0 w = \alpha$ on a short time interval.
In Section \ref{sec_moving_boundary_condition}, we show that the solution $s$ is strictly decreasing for a short time under assumption (C), and prove \red{Theorems} \ref{thm_main1}--\ref{thm_main2}.

\section{Mild solution} \label{sec_mild_solution}
\subsection{Derivation and formulas}
We begin by presenting the derivation of \eqref{eq_mild_solution_w}.
\begin{proposition}
    Assume \eqref{eq_assumptions_w0_f_in_def_of_mild_solution}.
    Let
    \begin{align*}
        (w, s)
        \in C([0, T); BUC(\Omega))
        \cap C((0, T); BUC^1(\Omega))
        \times C[0, T)
    \end{align*}
    be the solution to \eqref{eq_w}.
    Then $w$ satisfies \eqref{eq_mild_solution_w}\red{.}
\end{proposition}
\begin{proof}
    Applying the Laplace transform $\mathcal{L}$ with respect to $y$ to both sides of \eqref{eq_w}, we obtain
    \begin{align*}
        \partial_t \hat{w}
        - f(s, t)
        + \lambda \alpha
        - \lambda^2 \hat{w}
        - \lambda \dot{s} \hat{w}
        = e^{\lambda s(t)} \widehat{R_{s(t)}\partial_x f}(\lambda, t),
    \end{align*}
    where
    \begin{align*}
        \hat{w} (\lambda, t)
        = \mathcal{L} w, \quad
        R_\xi f
        =  \chi_{(\xi, \infty)} f.
        %R_\xi f
        %= \chi_{(0, \xi)} f, \quad
        %R_\xi f
        %=  (I - R_\xi) f.
    \end{align*}
    Solving the ODE, we deduce that
    \begin{align*}
        \partial_t (e^{- \lambda^2 t - \lambda s(t)} \hat{w}(t))
        & = e^{- \lambda^2 t - \lambda s(t)} f(s(t), t)
        - e^{- \lambda^2 t - \lambda s(t)} \lambda \alpha\\
        & + e^{- \lambda^2 t} \widehat{{R}_{s(t)} \partial_x f}(\lambda, t)
    \end{align*}
    and then
    \begin{align*}
        \hat{w}(t)
        & = e^{\lambda^2 t + \lambda s(t)} \hat{w_0}
        + \int_0^t
        e^{\lambda^2 (t - \tau) + \lambda (s(t) - s(\tau))} f(s(\tau), \tau)
        d\tau\\
        & - \int_0^t
        e^{\lambda^2 (t - \tau) + \lambda (s(t) - s(\tau))} \lambda \alpha
        d\tau
        + \int_0^t
            e^{\lambda^2 (t - \tau) + \lambda s(t)} \widehat{R_{s(\tau)} \partial_x f}(\lambda, \tau)
        d\tau.
    \end{align*}
    Since the support of $w_0$ is $[0, \infty)$, we can identify it with its zero extension, and thus obtain the solution formula \eqref{eq_mild_solution_w} by applying the inverse Laplace transform.
\end{proof}

\begin{remark}
    The reader might expect that the kernel contains a correction term arising from reflection, as in the case of the half-space, but this is not the case.
    Such a correction term is usually derived to cancel the integral term involving the boundary value of $\partial_y w$.
    However, in our problem, the boundary value of $\partial_y w$ is determined by $\dot{s}$ and $f(s(t),t)$, which follows from the Dirichlet boundary condition for $u$.
    Therefore, the reflection term is replaced by these terms.
\end{remark}
Let us split the solution formula as
\begin{align} \label{eq_def_of_w}
    \begin{split}
        w(y, t)
        & = \int_0^\infty
            G_t (y - z + s(t)) w_0(z)
        dz\\
        & + \int_0^t
            G_{t - \tau} (y + s(t) - s(\tau)) f(s(\tau), \tau)
        d\tau\\
        & - \alpha \int_0^t
            \partial_y G_{t - \tau}(y + s(t) - s(\tau))
        d\tau\\
        & + \int_0^{t}
            \int_0^\infty
                G_{t - \tau} (y - z + s(t) - s(\tau))
                \partial_z f(z + s(\tau), \tau)
            dz
        d\tau\\
        & =: w_1(s; y, t) + w_2(s; y, t) + w_3(s; y, t) + w_4(s; y, t).
    \end{split}
\end{align}
The terms $w_1$ and $w_3$ contain boundary layer effects.
The remaining terms can be interpreted as forcing terms.
To simplify notation, we write $w_j(s;t)$ for $w_j(s;y,t)$ when the dependence on $s$ needs to be specified.
When $s$ is clear from the context, we simply write $w_j(y,t)$ or $w_j(t)$ instead of $w_j(s;y,t)$.

Using \red{the} identities
\begin{gather} \label{eq_formula_of_kernel}
    \begin{split}
        G_t(y + s(t))
        = \int_0^\infty
            \partial_y G_t
            (
                y
                - z
                + s(t)
            )
        dz,\\
        %& = - \int_0^\infty
        %    \partial_z G_t
        %    (
        %        y
        %        - z
        %        + s(t)
        %    )
        %dz
        %- G(y + s(t))
        %= 0,\\
        \left(
            \partial_\tau
            + \partial_y^2
            + \dot{s}(\tau) \partial_y
        \right)
        G_{t - \tau}(y + s(t) - s(\tau))
        = 0,
    \end{split}
\end{gather}
for $y, t > 0$, we observe that
\begin{align} \label{eq_formula_w_1_w_3}
    \begin{split}
        & \gamma_0 \partial_y
        \left(
              w_1(s; t)
            + w_3(s; t)
        \right)\\
        %& = \gamma_0 \int_0^\infty
        %    \partial_y G_t (y - z + s(t)) w_0(z)
        %dz
        %- \gamma_0 \alpha \int_0^t
        %    \partial_y^2 G_{t - \tau}(y + s(t) - s(\tau))
        %d\tau\\
        & = \gamma_0 \int_0^\infty
            \partial_y G_t (y - z + s(t))
            \left(
                w_0(z)
                - w_0(y)
            \right)
        dz
        + \gamma_0 G_t(y + s(t)) w_0(y)\\
        & - \gamma_0 \alpha G_{t}(y + s(t))
        + \gamma_0 \alpha \int_0^t
            \dot{s}(\tau) \partial_y G_{t - \tau}(y + s(t) - s(\tau))
        d\tau\\
        %& = \int_0^\infty
        %    \partial_y G_t (- z + s(t))
        %    \left(
        %        w_0(z)
        %        - \alpha
        %    \right)
        %dz
        %+ G_t(s(t)) \alpha\\
        %& - \alpha G_{t}(s(t))
        %+ \alpha \int_0^t
        %    \dot{s}(\tau) \partial_y G_{t - \tau}(s(t) - s(\tau))
        %d\tau\\
        & = \int_0^\infty
            \partial_y G_t (- z + s(t))
            \left(
                w_0(z)
                - \alpha
            \right)
        dz\\
        & + \alpha \int_0^t
            \dot{s}(\tau) \partial_y G_{t - \tau}(s(t) - s(\tau))
        d\tau.
    \end{split}
\end{align}
We find from \eqref{eq_formula_w_1_w_3} that
\begin{align} \label{eq_formula_diff_w_1_w_3}
    \begin{split}
    & \gamma_0 \partial_y
        \left(
              w_1(s_2; t)
            - w_1(s_1; t)
            + w_3(s_2; t)
            - w_3(s_1; t)
        \right)\\
        & = \int_0^\infty
            \partial_y
            \left(
                  G_t (- z + s_2(t))
                - G_t (- z + s_1(t))
            \right)
            \left(
                w_0(z)
                - \alpha
            \right)
        dz\\
        & + \alpha \int_0^t
            (
                  \dot{s}_2(\tau)
                - \dot{s}_1(\tau)
            )
            \partial_y G_{t - \tau}(s_2(t) - s_2(\tau))
        d\tau\\
        & + \alpha \int_0^t
            \dot{s}_1(\tau) \partial_y
            \left(
                  G_{t - \tau}(s_2(t) - s_2(\tau))
                - G_{t - \tau}(s_1(t) - s_1(\tau))
            \right)
        d\tau\\
        & = \sigma(t)
        \int_0^\infty
            \int_0^1
                \partial_y^2 G_t (- z + s_\mu(t))
            d\mu
            \left(
                w_0(z)
                - \alpha
            \right)
        dz\\
        &
            + \alpha \int_0^t
            \dot{\sigma}(\tau)
            \partial_y G_{t - \tau}
            (
                s_2(t)
                - s_2(\tau)
            )
        d\tau\\
        &
            + \alpha \int_0^t
            \dot{s}_1(\tau)
            (
                \sigma(t)
                - \sigma(\tau)
            )
            \int_0^1
                \partial_y^2 G_{t - \tau}
                (
                    s_\mu(t)
                    - s_\mu(\tau)
                )
            d\mu
        d\tau,
    \end{split}
\end{align}
where
\begin{align*}
    \sigma
    &:= s_2
    - s_1,\\
    s_\mu
    &:= \mu s_2
    + (1 - \mu) s_1.
\end{align*}
If $w_0$ is differentiable\red{, we see} that
%\begin{align*}
%    \partial_y (w_1 + w_3)
%    = \alpha G_t(y + s(t))
%    + \int_0^\infty
%        G_t (y - z + s(t)) \partial_z w_0(z)
%    dz
%\end{align*}
%\begin{align*}
%    \partial_y w_3
%    = - \alpha G_t(y + s(t) - s(\tau))
%    + \int_0^t
%        \dot{s}(\tau) \partial_y G_t (y + s(t) - s(\tau))
%    d\tau
%\end{align*}
\begin{align} \label{eq_formula_w_1_w_3_when_w0_is_differentiable}
    \begin{split}
        & \partial_y
        \left(
              w_1(s; t)
            + w_3(s; t)
        \right)\\
        & = \alpha G_t(y + s(t))
        + \int_0^\infty
            G_t (y - z + s(t)) \partial_z w_0(z)
        dz\\
        & - \alpha G_t(y + s(t))
        + \alpha \int_0^t
            \dot{s}(\tau) \partial_y G_{t - \tau}(y + s(t) - s(\tau))
        d\tau\\
        & = \int_0^\infty
            G_t (y - z + s(t)) \partial_z w_0(z)
        dz\\
        & + \alpha \int_0^t
            \dot{s}(\tau) \partial_y G_{t - \tau}(y + s(t) - s(\tau))
        d\tau.
    \end{split}
\end{align}
\begin{remark} \label{rmk_vanishing_boundary_layer}
    In the cancellation of $\alpha G_t(y+s(t))$ in \eqref{eq_formula_w_1_w_3} and \eqref{eq_formula_w_1_w_3_when_w0_is_differentiable}, the boundary condition $\gamma_0 w_0 = \alpha$ is essential.
    Moreover, in order for $s^\prime$ to make sense up to $t=0$, the term $\alpha G_t(y+s(t))$ must vanish, since this term behaves like $O(t^{-\frac{1}{2}})$.
    Therefore, some differentiability of $w_0$ is essentially necessary for continuity at $t=0$.
    This observation is related to Theorem \ref{thm_main2}.
\end{remark}
Similarly to \eqref{eq_formula_diff_w_1_w_3}, we also  have
\begin{align} \label{eq_formula_diff_w_1_w_3_when_w0_is_differentiable}
    \begin{split}
    & \gamma_0 \partial_y
        \left(
              w_1(s_2; t)
            - w_1(s_1; t)
            + w_3(s_2; t)
            - w_3(s_1; t)
        \right)\\
        & = \sigma(t)
        \int_0^\infty
            \int_0^1
                \partial_y G_t (- z + s_\mu(t))
            d\mu
            \partial_y w_0(z)
        dz\\
        & + \alpha \int_0^t
            \dot{\sigma}(\tau)
            \partial_y G_{t - \tau}
            (
                s_2(t)
                - s_2(\tau)
            )
        d\tau\\
        & 
        + \alpha \int_0^t
            \dot{s}_1(\tau)
            (
                \sigma(t)
                - \sigma(\tau)
            )
            \int_0^1
                \partial_y^2 G_{t - \tau}
                (
                    s_\mu(t)
                    - s_\mu(\tau)
                )
            d\mu
        d\tau.
    \end{split}
\end{align}

We next consider estimates for $w_j$.
To estimate the kernels, we observe
\begin{align*}
    & \vert
        s(t_2)
        - s(t_1)
    \vert\\
    & \leq \int_0^1
        \vert
            \dot{s}(\tau t_2 + (1 - \tau) t_1)
        \vert
    d\tau
    (
        t_2
        - t_1
    )
    \leq \frac{1}{\theta}
    (
        t_2^{\theta}
        - t_1^{\theta}
    )
    \Vert
        s
    \Vert_{X^\theta(T)}\\
    & \leq \frac{1}{\theta}
    \min
    \left(
        t_1^{- 1 + \theta}
        (
            t_2
            - t_1
        ),
        (
            t_2
            - t_1
        )^\theta
    \right)
    \Vert
        s
    \Vert_{X^\theta(T)}
\end{align*}
for all $s \in X^\theta(T)$ and $0 < t_1 < t_2$.
Since $s(0) = 0$, we also have
\begin{align*}
    \vert
        s(t)
    \vert
    = \left \vert
        \int_0^t
            \dot{s}(\tau)
        d\tau
    \right \vert
    \leq \frac{t^\theta}{\theta}
    \Vert
        s
    \Vert_{X^\theta(T)}.
\end{align*}
We estimate the following integral that typically appears in convolution estimates.
A typical case is $\varphi = 1$ and $r = \infty$ ($r^\prime = 1$).
\begin{proposition} \label{prop_beta_type}
    Let $r \in [1, \infty]$, $a, b, c \geq 0$ satisfy $a - b > - \frac{1}{r^\prime}$ and $c < \frac{1}{r^\prime}$.
    Let $s \in X^\theta(T)$.
    Then
    \begin{align*}
        \int_0^t
            \vert
                s(t)
                - s(\tau)
            \vert^a
            (
                t
                - \tau
            )^{-b}
            \tau^{-c}
            \varphi(\tau)
        d\tau
        \leq C t^{1 - \frac{1}{r} + a \theta - b - c}
        \Vert
            \varphi
        \Vert_{L^r(0, t)}.
    \end{align*}
    for all $\varphi \in L^r(0, T)$.
    A simple change of variables yields
    \begin{align*}
        \int_0^t
            \vert
                s(t)
                - s(t - \tau)
            \vert^a
            \tau^{-b}
            (t - \tau)^{-c}
            \varphi(t - \tau)
        d\tau
        \leq C t^{1 - \frac{1}{r} + a \theta - b - c}
        \Vert
            \varphi
        \Vert_{L^r(0, t)}.
    \end{align*}
\end{proposition}
\begin{proof}
    Let $\overline{a} = a r^\prime$, $\overline{b} = b r^\prime$, and $\overline{c} = c r^\prime$.
    We deduce that
    \begin{align*}
        & \left(
            \int_0^\frac{t}{2}
            + \int_\frac{t}{2}^t
        \right)
            \vert
                s(t)
                - s(\tau)
            \vert^{\overline{a}}
            (
                t
                - \tau
            )^{- \overline{b}}
            \tau^{- \overline{c}}
        d\tau\\
        & \leq C
        \Vert
            s
        \Vert_{X^\theta(T)}
        \left(
            \int_0^\frac{t}{2}
                (t - \tau)^{\overline{a} \theta - \overline{b}} \tau^{- \overline{c}}
            d\tau
            + \int_\frac{t}{2}^t
                (t - \tau)^{\overline{a} - \overline{b}} \tau^{- \overline{a} + \theta \overline{a} - \overline{c}}
            d\tau
        \right)\\
        %& \leq C t^{1 + \overline{a} \theta - \overline{b} - \overline{c}}
        %+ C t^{1 + \overline{a} - \overline{b} - \overline{a} + \theta a - \overline{c}}\\
        & \leq C
        \Vert
            s
        \Vert_{X^\theta(T)}
        t^{1 + \overline{a} \theta - \overline{b} - \overline{c}}.
    \end{align*}
    The estimate follows from the H\"{o}lder inequality.
\end{proof}
In the following subsections, we derive estimates for $w_1$, $w_3$, $w_2$, and $w_4$.
To estimate the trace $\gamma_0 \partial_y w$, we consider $w_1$ and $w_3$ together.

\subsection{Estimates for $w_1$ and $w_3$} \label{subsec_w1_w3}
We first consider the case $s \in X^{\frac12}(T)$ and then the case $s \in Y(T)$.
The $L^q$-estimates are derived separately.
\begin{proposition} \label{prop_estimates_for_w1}
    Let $q \in [1, \infty]$.
    Let $w_0 \in L^q(\Omega)$ and $s \in C[0, \infty)$.
    Then
    \begin{align*}
        \sup_{0 < t < T}
            \Vert
                w_1(t)
            \Vert_{L^q(\Omega)}
        + \sup_{0 < t < T}
            t^{\frac{1}{2}}
            \Vert
                \partial_y w_1(t)
            \Vert_{L^q(\Omega)}
        \leq C \Vert
            w_0
        \Vert_{L^{q}(\Omega)}.
    \end{align*}
\end{proposition}
\begin{proof}
    %Continuity is a consequence of the Lebesgue dominated convergence theorem.
    The Young inequality yields the estimate.
\end{proof}
\begin{proposition} \label{prop_estimates_for_w3}
    Let $q \in [1, \infty]$ and $s \in X^\frac{1}{2}(T)$.
    Then
    \begin{align*}
        \sup_{0 < t < T}
            \Vert
                w_3(t)
            \Vert_{L^q(\Omega)}
        & \leq C \alpha
        \left(
            1
            + \Vert
                s
            \Vert_{X^\frac{1}{2}(T)}
        \right),\\
        \sup_{0 < t < T}
            %t^{\frac{1}{2} - \frac{1}{2q}}
            t^{\frac{1}{2}}
            \Vert
                \partial_y w_3(t)
            \Vert_{L^q(\Omega)}
        & \leq C \alpha
        \left(
            1
            +
            \Vert
                s
            \Vert_{X^\frac{1}{2}(T)}
            +
            \Vert
                s
            \Vert_{X^\frac{1}{2}(T)}^2
        \right).
    \end{align*}
\end{proposition}
\begin{proof}
    We first observe that
    \begin{align*}
        w_3(t)
        & = 
        - \alpha \int_0^t
            \partial_y G_\tau (
                y
            )
        d\tau\\
        & - \alpha \int_0^t
            \int_0^1
                \partial_y^2 G_\tau (
                    y
                    - \mu
                    (
                        s(t)
                        - s(t - \tau)
                    )
                )
            d\mu
            (
                s(t)
                - s(t - \tau)
            )
        d\tau\\
        & =: I_1 + I_2.
    \end{align*}
    By a change of variables, we see that
    \begin{align} \label{eq_int_dG_tau_dtau_I1}
        \begin{split}
            \Vert
                I_1
            \Vert_{L^q(\Omega)}
            & \leq \alpha
            \left \Vert
                \int_0^t
                    \vert
                        \partial_y G_\tau (
                            y
                        )
                    \vert
                d\tau
            \right \Vert_{L^q(\Omega)}\\
            & = \alpha
            \left \Vert
                \frac{1}{2 \sqrt{\pi}}
                \int_{\frac{y}{\sqrt{t}}}^\infty
                    e^{- \frac{\mu^2}{4}}
                d\mu
            \right \Vert_{L^q(\Omega)}
            \leq C \alpha t^{\frac{1}{2q}}
        \end{split}
    \end{align}
    for all $q \in [1, \infty]$.
    Since $s \in X^\frac{1}{2}(T)$, we use Proposition \ref{prop_beta_type} to deduce that
    \begin{align} \label{eq_int_dG_tau_dtau_I2}
        \begin{split}
            \Vert
                I_2
            \Vert_{L^q(\Omega)}
            & \leq C \alpha
            \int_0^t
                \tau^{-1 - \frac{1}{2} \left(1 - \frac{1}{q}\right)}
                \vert
                      s(t)
                    - s(t - \tau)
                \vert
            d\tau\\
            & \leq C \alpha
            \Vert
                s
            \Vert_{X^\frac{1}{2}(T)}
            t^{\frac{1}{2q}}
        \end{split}
    \end{align}
    for all $q \in [1, \infty]$.

    We find from the second identity of \eqref{eq_formula_of_kernel} that
    \begin{align*}
        & \partial_y \int_0^t
            \partial_y G_{t - \tau}(y + s(t) - s(\tau))
        d\tau\\
        & = G_{t} (y + s(t))
        - \int_0^t
            \dot{s}(\tau)
            \partial_y G_{t - \tau} (y + s(t) - s(\tau))
        d\tau.
    \end{align*}
    This formula implies the $L^1$-estimate
    \begin{align*}
        \left \Vert
            \partial_y \int_0^t
                \partial_y G_{t - \tau}(y + s(t) - s(\tau))
            d\tau
        \right \Vert_{L^1(\Omega)}
        \leq C
        (
            1
            +
            \Vert
                s
            \Vert_{X^\frac{1}{2}(T)}
        ).
    \end{align*}
    We next derive the $L^\infty$-estimate.
    We continue the estimate via a change of variables as
    \begin{align} \label{eq_bound_for_int_sdot_dy_G_tau}
        \begin{split}
            & \left \vert
                \int_0^t
                    \dot{s}(\tau)
                    \partial_y G_{t - \tau} (y + s(t) - s(\tau))
                d\tau
            \right \vert\\
            %& = \int_0^t
            %    \dot{s}(t - \tau)
            %    \partial_y G_{\tau} (y + s(t) - s(t - \tau))
            %d\tau\\
            & \leq \left \vert
                \int_0^t
                    \dot{s}(t - \tau)
                    \partial_y G_{\tau} (y)
                d\tau
            \right \vert\\
            & + \left \vert
                \int_0^t
                    \int_0^1
                        \dot{s}(t - \tau)
                            \partial_y^2 G_{\tau} (y + \mu (s(t) - s(t - \tau)))
                    d\mu
                    (
                        s(t)
                        - s(t - \tau)
                    )
                d\tau
            \right \vert.
        \end{split}
    \end{align}
    The first term is bounded by the change of variables $\mu = \frac{y}{\sqrt{\tau}}$ ($\frac{d\tau}{\tau} = - \frac{2 d\mu}{\mu}$) as
    \begin{align*}
        & \left \vert
            \int_0^t
                \dot{s}(t - \tau)
                \partial_y G_{\tau} (y)
            d\tau
        \right \vert
        \leq \frac{1}{2 \sqrt{\pi}}
        \Vert
            s
        \Vert_{X^\frac{1}{2}(T)}
        \left \vert
            \int_0^t
                (t - \tau)^{- \frac{1}{2}}
                \frac{y}{\tau^{\frac{3}{2}}}
                e^{- \frac{y^2}{4\tau}}
            d\tau
        \right \vert\\
        & = \frac{1}{\sqrt{\pi}}
        \Vert
            s
        \Vert_{X^\frac{1}{2}(T)}
        \int_{\frac{y}{\sqrt{t}}}^\infty
            \left(
                t
                - \frac{y^2}{\mu^2}
            \right)^{- \frac{1}{2}}
            \mu
            e^{- \frac{\mu^2}{4}}
        \frac{d\mu}{\mu}\\
        %= \frac{1}{\sqrt{\pi}}
        %\Vert
        %    s
        %\Vert_{X^\frac{1}{2}(T)}
        %\int_{\frac{y}{\sqrt{t}}}^\infty
        %    \left(
        %        t
        %        - \frac{y^2}{\mu^2}
        %    \right)^{- \frac{1}{2}}
        %    e^{- \frac{\mu^2}{4}}
        %d\mu\\
        & \leq
        \frac{
            \Vert
                s
            \Vert_{X^\frac{1}{2}(T)}
        }{
            \sqrt{\pi t}
        }
        \int_{u}^\infty
            \frac{
                \mu e^{- \frac{\mu^2}{4}}
            }{
                \sqrt{\mu^2 - u^2}
            }
        d\mu
        \Biggr \vert_{u = \frac{y}{\sqrt{t}}}.
    \end{align*}
    Changing variables again by setting $\nu = \sqrt{\mu^2 - u^2}$ ($d\nu = \mu d\mu / \sqrt{\mu^2 - u^2}$), we deduce that
    \begin{align*}
        \left \vert
            \int_0^t
                \dot{s}(t - \tau)
                \partial_y G_{\tau} (y)
            d\tau
        \right \vert
        & \leq
        \frac{
            \Vert
                s
            \Vert_{X^\frac{1}{2}(T)}
        }{
            \sqrt{\pi t}
        }
        \int_{0}^\infty
            e^{- \frac{\nu^2 + u^2}{4}}
        d\nu
        \Biggr \vert_{u = \frac{y}{\sqrt{t}}}\\
        & = \frac{
            \Vert
                s
            \Vert_{X^\frac{1}{2}(T)}
        }{
            \sqrt{t}
        }
        e^{- \frac{y^2}{4t}}.
    \end{align*}
    We continue to estimate the second term.
    Since $r e^{-r}$ is bounded for $r > 0$, using Proposition \ref{prop_beta_type}, we have
    \begin{align} \label{eq_int_sdot_d2Gdtau}
        \begin{split}
            & \left \vert
                \int_0^t
                    \int_0^1
                        \dot{s}(t - \tau)
                            \partial_y^2 G_{\tau} (y + \mu (s(t) - s(t - \tau)))
                    d\mu
                    (
                        s(t)
                        - s(t - \tau)
                    )
                d\tau
            \right \vert\\
            & \leq
            \Vert
                s
            \Vert_{X^\frac{1}{2}(T)}
            \int_0^1
                \int_0^t
                    (t - \tau)^{- \frac{1}{2}}
                    \left(
                        \frac{1}{\sqrt{\tau}^3}
                        + \frac{
                            \vert
                                y
                                + \mu
                                ( 
                                    s(t)
                                    - s(t - \tau)
                                )
                            \vert^2
                        }{
                            \sqrt{\tau}^5
                        }
                    \right)\\
            & \quad\quad\quad\quad\quad\quad\quad\quad\quad
                    \times e^{
                        - \frac{
                            \vert
                                y
                                + \mu
                                ( 
                                    s(t)
                                    - s(t - \tau)
                                )
                            \vert^2
                        }{
                            4 \tau
                        }
                    }
                    \vert
                        s(t)
                        - s(t - \tau)
                    \vert
                d\tau
            d\mu\\
            & \leq \frac{
                C
                \Vert
                    s
                \Vert_{X^\frac{1}{2}(T)}^2
            }{
                \sqrt{t}
            }.
        \end{split}
    \end{align}
    By interpolation, we obtain the desired estimate.
\end{proof}
The estimates for the traces of the derivatives must be treated together; see Remark \ref{rmk_vanishing_boundary_layer}.
\begin{proposition} \label{prop_estimates_for_w1_w3}
    Let $w_0 \in BUC(\Omega)$ satisfy $\gamma_0 w_0 = \alpha$ and $s \in X^{\frac{1}{2}}(T)$.
    Then
    \begin{align*}
        \sup_{0 < t < T}
            t^{\frac{1}{2}}
            \vert
                \gamma_0 \partial_y
                (
                      w_1(s; t)
                    + w_3(s; t)
                )
            \vert
        \leq
        C_{w_0, T}
        + C \alpha
            \Vert
            s
        \Vert_{X^\frac{1}{2}(T)}^2
    \end{align*}
    for some positive constant $C_{w_0, T} \leq \frac{\Vert w_0 - \alpha \Vert_\infty}{\sqrt{\pi}}$ satisfying $C_{w_0, T} = o(1)$.
\end{proposition}
\begin{proof}
    The formula \eqref{eq_formula_w_1_w_3} and Proposition \ref{prop_beta_type} lead to
    \begin{align} \label{eq_estimate_gamma0_dw1w_3}
        \begin{split}
            & \vert
                \gamma_0 \partial_y
                \left(
                      w_1(s; t)
                    + w_3(s; t)
                \right)
            \vert\\
            & \leq \int_0^\infty
                \vert
                    \partial_y G_t (- z + s(t))
                    \left(
                        w_0(z)
                        - \alpha
                    \right)
                \vert
            dz\\
            & + C \alpha
            \int_0^t
                \vert
                    \dot{s}(\tau)
                \vert
                \vert
                    s(t)
                    - s(\tau)
                \vert
                (
                    t
                    - \tau
                )^{- \frac{3}{2}}
            d\tau\\
            & \leq t^{- \frac{1}{2}}
            \left(
                \frac{
                    \Vert
                        w_0
                        - \alpha
                    \Vert_\infty
                }{
                    \sqrt{\pi}
                }
                + C \alpha
                    \Vert
                    s
                \Vert_{X^\frac{1}{2}(T)}^2
            \right).
        \end{split}
    \end{align}
    For small $\varepsilon > 0$, there exists \red{a} small $\delta > 0$ \red{such that} if $0 < z < \delta$, \red{then} $\vert w_0(z) - \alpha \vert < \varepsilon$.
    For $z > \delta$ and sufficiently small $t > 0$, $\vert s(t) \vert < \frac{\delta}{2}$ and $\vert z + s(t) \vert > \frac{\delta}{2}$.
    The integrand $\vert t^\frac{1}{2} \partial_y G_t (- z + s(t)) \vert$ is bounded and exponentially decaying at \red{infinity} with order at least $\sup_{0 < \delta < 1} \vert \partial_y G_1 (- z + \frac{\delta}{2}) \vert$ and tends to zero as $t \rightarrow 0$.
    By Lebesgue's dominated convergence theorem, we see that
    \begin{align*}
        & \left(
            \int_0^\delta
            + \int_\delta^\infty
        \right)
            \vert
                t^\frac{1}{2}
                \partial_y G_t (- z + s(t))
            \vert
            \vert
                \left(
                    w_0(z)
                    - \alpha
                \right)
            \vert
        dz\\
        & 
        \leq
        \varepsilon
        \int_\Real
            \vert
                t^\frac{1}{2}
                \partial_y G_t (z)
            \vert
        dz
        + C
        \int_\delta^\infty
            \vert
                t^\frac{1}{2}
                \partial_y G_t (- z + s(t))
            \vert
        dz
        \\
        &\leq C \varepsilon.
    \end{align*}
\end{proof}

\begin{proposition} \label{prop_estimates_for_diff_w1_w_3}
    Let $w_0 \in BUC(\Omega)$ satisfy $\gamma_0 w_0 = \alpha$.
    Let $s_1, s_2 \in X^{\frac{1}{2}}(T)$.
    Then
    \begin{align*}
        & \sup_{0 < t < T}
            t^{\frac{1}{2}}
            \vert
                \gamma_0 \partial_y
                (
                      w_1(s_2; t)
                    - w_1(s_1; t)
                    + w_3(s_2; t)
                    - w_3(s_1; t)
                )
            \vert\\
        & \leq \left(
            C_{w_0, T}^\prime
            + C \alpha
            \Vert
                s_2
            \Vert_{X^\frac{1}{2}(T)}
        \right)
        \Vert
            \sigma
        \Vert_{X^\frac{1}{2}(T)}
    \end{align*}
    for some constant $C_{w_0, T}^\prime = o(1)$.
\end{proposition}
\begin{proof}
    We deduce from \eqref{eq_formula_diff_w_1_w_3} and the spatially uniform estimate for the derivatives of the Gaussian kernel that
    \begin{align*}
        & \vert
            \gamma_0 \partial_y
            (
                  w_1(s_2; t)
                - w_1(s_1; t)
                + w_3(s_2; t)
                - w_3(s_1; t)
            )
        \vert\\
        & \leq
        \vert
            \sigma(t)
        \vert
        \int_0^\infty
            \int_0^1
                \vert
                    \partial_y^2 G_t (- z + s_\mu(t))
                \vert
            d\mu
            \vert
                w_0(z)
                - \alpha
            \vert
        dz\\
        & + \alpha
        \int_0^t
            \vert
                \dot{\sigma}(\tau)
            \vert
            \vert
                \partial_y G_{t - \tau}
                (
                    s_2(t)
                    - s_2(\tau)
                )
            \vert
        d\tau\\
        & + \alpha \int_0^t
            \vert
                \dot{s}_1(\tau)
            \vert
            \vert
                \sigma(t)
                - \sigma(\tau)
            \vert
            \int_0^1
                \vert
                    \partial_y^2 G_{t - \tau}
                    (
                        s_\mu(t)
                        - s_\mu(\tau)
                    )
                \vert
            d\mu
        d\tau\\
        & \leq
        C_1 t^{- \frac{1}{2}}
        \Vert
            \sigma
        \Vert_{X^\frac{1}{2}(T)}
        \\
        &
        + C_2 \alpha
        \Vert
            \sigma
        \Vert_{X^\frac{1}{2}(T)}
        \int_0^t
            \tau^{- \frac{1}{2}}
            (
                t
                - \tau
            )^{- \frac{1}{2}}
        d\tau
        \\
        &
        + C_3 \alpha
        \Vert
            s_1
        \Vert_{X^\frac{1}{2}(T)}
        \int_0^t
            \tau^{- \frac{1}{2}}
            \vert
                \sigma(t)
                - \sigma(\tau)
            \vert
            (
                t
                - \tau
            )^{- \frac{3}{2}}
        d\tau
        \\
        &
        \leq \frac{1}{\sqrt{t}}
        \left(
            C_1
            + \alpha C_2 t^\frac{1}{2}
            %\frac{2 \sqrt{2}}{\sqrt{t \pi e}}
            %\Vert
            %    w_0
            %    - \alpha
            %\Vert_\infty
            + C \alpha
            \Vert
                s_1
            \Vert_{X^\frac{1}{2}(T)}
        \right)
        \Vert
            \sigma
        \Vert_{X^\frac{1}{2}(T)}.
    \end{align*}
    \red{As in} the proof of Proposition \ref{prop_estimates_for_w1_w3}, we find \red{that} $C_1 = o(1)$ and \red{that $C_1$ is} bounded by $\frac{2 \sqrt{2}}{\sqrt{\pi e}} \Vert w_0 - \alpha \Vert_\infty$.
\end{proof}

We consider the case \red{where} $w_0$ is differentiable and $s \in Y(T)$\red{.}
\begin{proposition} \label{eq_estimate_for_w1_w3_if_w0_is_differentiable}
    Assume that $w_0 \in BUC^1(\Omega)$ and $s \in Y(T)$.
    Then
    \begin{align*}
        \sup_{0 < t < T}
            \Vert
                \partial_y(
                    w_1(t)
                    + w_3(t)
                )
            \Vert_\infty
        \leq C
        \left(
            \Vert
                \partial_y w_0
            \Vert_{\infty}
            + \alpha
            (
                \Vert
                    s
                \Vert_{Y(T)}
                + \Vert
                    s
                \Vert_{Y(T)}^2
            )
        \right)
    \end{align*}
    and
    \begin{align*}
        \sup_{0 < t < T}
            \vert
                \gamma_0\partial_y
                 (
                    w_1(t)
                    + w_3(t)
                )
            \vert
        \leq C
        \left(
            \Vert
                \partial_y w_0
            \Vert_{\infty}
            + \alpha  T^\frac{1}{2}
            \Vert
                s
            \Vert_{Y(T)}^2
        \right).
    \end{align*}
\end{proposition}
\begin{proof}
    We deduce from the formula \eqref{eq_formula_w_1_w_3_when_w0_is_differentiable} and the Young inequality that
    \begin{align*}
        & \vert
            \partial_y(
                w_1(t)
                + w_3(t)
            )
        \vert
        \leq C \Vert
            \partial_y w_0
        \Vert_{\infty}
        + \alpha
        \left \vert
            \int_0^t
                \dot{s}(\tau) \partial_y G_{t - \tau}(y + s(t) - s(\tau))
            d\tau
        \right \vert.
    \end{align*}
    The second term is bounded in the same way as \eqref{eq_bound_for_int_sdot_dy_G_tau}, \red{and} then we obtain the first estimate.
    The second estimate is derived in the same way as \eqref{eq_estimate_gamma0_dw1w_3}.
\end{proof}
\begin{proposition} \label{eq_estimate_for_diff_w1_w3_if_w0_is_differentiable}
    Let $w_0 \in BUC^1(\Omega)$ and $s_1, s_2 \in Y(T)$.
    %Assume that $T < 1$.
    Then
    \begin{align*}
        & \sup_{0 < t < T}
            \vert
                \gamma_0 \partial_y(
                      w_1(s_2; t)
                    + w_3(s_2; t)
                    - w_1(s_1; t)
                    - w_3(s_1; t)
                )
            \vert\\
        & \leq C T^\frac{1}{2}
        \left(
            \Vert
                \partial_y w_0
            \Vert_{BUC(\Omega)}
            + \alpha
            \sum_{j=1,2}
                \Vert
                    s_j
                \Vert_{Y(T)}
        \right)
        \Vert
            \sigma
        \Vert_{Y(T)}.
    \end{align*}
\end{proposition}
\begin{proof}
    We deduce from \eqref{eq_formula_diff_w_1_w_3_when_w0_is_differentiable} and the spatially uniform estimate for the derivatives of the Gaussian kernel that
    \begin{align*}
        & \vert
            \gamma_0 \partial_y
            \left(
                  w_1(s_2; t)
                - w_1(s_1; t)
                + w_3(s_2; t)
                - w_3(s_1; t)
            \right)
        \vert\\
        & \leq C t^\frac{1}{2}
        \Vert
            \sigma
        \Vert_{Y(T)}
        \Vert
            \partial_y w_0
        \Vert_{BUC(\Omega)}
        + C \alpha
        \Vert
            \sigma
        \Vert_{Y(T)}
        \Vert
            s_2
        \Vert_{Y(T)}
        \int_0^t
            (
                t
                - \tau
            )^{- \frac{1}{2}}
        d\tau\\
        & + C \alpha
        \Vert
            s_1
        \Vert_{Y(T)}
        \Vert
            \sigma
        \Vert_{Y(T)}
        \int_0^t
            (
                t
                - \tau
            )^{- \frac{1}{2}}
        d\tau\\
        & \leq C T^\frac{1}{2}
        \left(
            \Vert
                \partial_y w_0
            \Vert_{BUC(\Omega)}
            + \alpha
            \sum_{j=1,2}
                \Vert
                    s_j
                \Vert_{Y(T)}
        \right)
        \Vert
            \sigma
        \Vert_{Y(T)}.
    \end{align*}
\end{proof}

\subsection{Estimates for $w_2$} \label{subsec_w2}
We next estimate the boundary forcing term $w_2$.
We begin by considering the case where $s \in X^\frac{1}{2}(T)$.
\begin{proposition} \label{prop_estimates_for_w2}
    Let $q \in [1, \infty]$, $r \in (2q, \infty]$, and $\eta \in [0,1 - \frac{1}{r})$ satisfy
    \begin{align*}
        \eta
        \leq \frac{1}{2} - \frac{1}{r} + \frac{1}{2 q}.
    \end{align*}
    If $f \in C(\Real \times (0, T))$ \red{satisfies}
    \begin{gather*}
        t^{\eta} f
        \in L^r(0, T; BUC(\Real)),
    \end{gather*}
    then
    \begin{gather*}
        \sup_{0 < t < T}
            \Vert
                w_2(t)
            \Vert_{L^{q}(\Omega)}
        + \sup_{0 < t < T}
            t^{\frac{1}{2}}
            \Vert
                \partial_y w_2(t)
            \Vert_{L^{q}(\Omega)}
        \leq C \Vert
            t^{\eta}f
        \Vert_{L^r(0, T; BUC(\Real))}
    \end{gather*}
    for $q \in [1, \infty)$.
    When $q = \infty$, \red{assume that} $f$ satisfies $C_f < \infty$.
    Then
    \begin{align*}
        & \sup_{0 < t < T}
            \Vert
                w_2(t)
            \Vert_\infty
        + \sup_{0 < t < T}
            t^{\frac{1}{2}}
            \Vert
                \partial_y w_2(t)
            \Vert_\infty\\
        & \leq C (
            C_f
            +
            \sup_{0 < t < T}
                t^\frac{1}{2}
                \Vert
                    f(t)
                \Vert_\infty
        )
        \left(
            1
            + \Vert
                s
            \Vert_{X^\frac{1}{2}(T)}
        \right).
    \end{align*}
\end{proposition}
\begin{proof}
    We find from the H\"{o}lder inequality and Proposition \ref{prop_beta_type} that
    \begin{align*}
        \Vert
            w_2(t)
        \Vert_{L^q(\Omega)}
        & \leq C \int_0^t
            (t - \tau)^{- \frac{1}{2} + \frac{1}{2q}}
            f(s(\tau), \tau)
        d\tau\\
        & \leq C t^{\frac{1}{2} - \frac{1}{r} + \frac{1}{2 q} - \eta}
        \Vert
            t^{\eta} f
        \Vert_{L^r(0, T; BUC(\Real))}
    \end{align*}
    for $q \in [1, \infty]$ and
    \begin{align*}
        \Vert
            \partial_y w_2(t)
        \Vert_{L^q(\Omega)}
        & \leq C t^{- \frac{1}{r} + \frac{1}{2 q} - \eta}
        \Vert
            t^{\eta} f
        \Vert_{L^r(0, T; BUC(\Real))}
    \end{align*}
    for $q \in [1, \infty)$.
    For the case that $q = \infty$, we split $\partial_y w_2$ by
    \begin{align} \label{eq_formula_dw_2}
        \begin{split}
            \partial_y w_2
            & = \int_0^t
                \partial_y G_{t - \tau} (y + s(t) - s(\tau))
                \left(
                    f(s(\tau), \tau)
                    - f(s(t), t)
                \right)
            d\tau\\
            & + f(s(t), t)
            \int_0^t
                \partial_y G_{t - \tau} (y + s(t) - s(\tau))
            d\tau.
        \end{split}
    \end{align}
    The estimate for the second \red{term} can be done in the same way as Proposition \ref{prop_estimates_for_w3}.
    The first term is bounded by
    \begin{align*}
        &\left \Vert
            \int_0^t
                \partial_y G_{t - \tau} (y + s(t) - s(\tau))
                \left(
                    f(s(\tau), \tau)
                    - f(s(t), t)
                \right)
            d\tau
        \right \Vert_{BUC(\Omega)}\\
        & \leq C C_f
        \left(
            1
            + \Vert
                s
            \Vert_{X^\frac{1}{2}(T)}
        \right)
        \int_0^t
            (t - \tau)^{- \frac{1}{2}} t^{- \frac{1}{2}} \tau^{- \frac{1}{2}}
        d\tau
        = O(t^{- \frac{1}{2}}).
    \end{align*}
    From these observations, we obtain the conclusion.
\end{proof}

\begin{proposition} \label{prop_estimates_for_diff_w2}
    Let $\varepsilon \in [0, \frac{1}{2}]$.
    Let $s_1, s_2 \in X^\frac{1}{2}(T)$ and
    \begin{align*}
        t^\frac{1}{2} f, \,
        t^{1 - \varepsilon} \partial_x f
        \in BUC(\Real \times (0, T)).
    \end{align*}
    Then
    \begin{align*}% \label{eq_estimate_w_2}
        \begin{split}
            & \sup_{0 < t < T}
                t^{\frac{1}{2}}
                \vert
                    \gamma_0 \partial_y
                    (
                          w_2(s_2; t)
                        - w_2(s_1; t)
                    )
                \vert\\
            & \leq C
            \Vert
                \sigma
            \Vert_{X^\frac{1}{2}(T)}
            \left(
                \Vert
                    t^\frac{1}{2} f
                \Vert_{BUC(\Real \times (0, T))}
                + T^\varepsilon
                \Vert
                    s_1
                \Vert_{X^\frac{1}{2}(T)}
                \Vert
                    t^{1 - \varepsilon} \partial_x f
                \Vert_{BUC(\Real \times (0, T))}
            \right).
        \end{split}
    \end{align*}
\end{proposition}
\begin{proof}
    We find from the formula
    \begin{align} \label{eq_formula_gamma_0_dw2s2_w2s1}
        \begin{split}
            & \gamma_0 \partial_y
                \left(
                    w_2(s_2; t)
                    - w_2(s_1; t)
                \right)\\
            & = \int_0^t
                (
                      \sigma(t)
                    - \sigma(\tau)
                )
                \int_0^1
                    \partial_y^2 G_{t - \tau} (s_\mu(t) - s_\mu(\tau))
                d\mu
                f(s_2(\tau), \tau)
            d\tau\\
            & + \int_0^t
                \partial_y G_{t - \tau} (s_1(t) - s_1(\tau))
                \sigma(\tau)
                \int_0^1
                    \partial_x f(s_\mu(\tau), \tau)
                d\mu
            d\tau
        \end{split}
    \end{align}
    and Proposition \ref{prop_beta_type} that
    \begin{align*}
        & \left \vert
            \gamma_0 \partial_y
            \left(
                w_2(s_2; t)
                - w_2(s_1; t)
            \right)
        \right \vert \\
        & \leq C
        \int_0^t
            \vert
                  \sigma(t)
                - \sigma(\tau)
            \vert
            \int_0^1
                \left[
                    \frac{
                        1
                    }{
                        (
                            t
                            - \tau
                        )^\frac{3}{2}
                    }
                    + \frac{
                        \vert
                            s_\mu(t)
                            - s_\mu(\tau)
                        \vert
                    }{
                        (
                            t
                            - \tau
                        )^\frac{5}{2}
                    }
                \right]
                e^{- \frac{
                    \vert
                        s_\mu(t)
                        - s_\mu(\tau)
                    \vert^2
                }{
                    4
                    (
                        t
                        - \tau
                    )
                }}
            d\mu
            \vert
                f(s_2(\tau), \tau)
            \vert
        d\tau\\
        & + C
        \int_0^t
            \frac{
                \vert
                    s_1(t)
                    - s_1(\tau)
                \vert
            }{
                (
                    t
                    - \tau
                )^\frac{3}{2}
            }
            \vert
                \sigma(\tau)
            \vert
            \int_0^1
                \vert
                    \partial_x f(s_\mu(\tau), \tau)
                \vert
            d\mu
        d\tau\\
        & \leq C t^{- \frac{1}{2}}
        \Vert
            \sigma
        \Vert_{X^\frac{1}{2}(T)}
        \Vert
            t^\frac{1}{2} f
        \Vert_{BUC(\Real \times (0, T))}\\
        & + C t^{- \frac{1}{2} + \varepsilon}
        \Vert
            \sigma
        \Vert_{X^\frac{1}{2}(T)}
        \Vert
            s_1
        \Vert_{X^\frac{1}{2}(T)}
        \Vert
            t^{1 - \varepsilon} \partial_x f
        \Vert_{BUC(\Real \times (0, T))}.
    \end{align*}
\end{proof}
We next consider the case where $s \in Y(T)$.
\begin{proposition} \label{prop_estimates_for_w2_for_smooth_f}
    Let $s \in Y(T)$, and $f \in Lip(\Real \times (0, T))$.
    %Assume that $T < 1$.
    Then
    \begin{gather*}
        \sup_{0 < t < T}
            \Vert
                \partial_y w_2(t)
            \Vert_{\infty}
        \leq C
        \Vert
            f
        \Vert_{BUC((0, T); BUC^1(\Real))}
        (
            1
            + T^\frac{1}{2}
            \Vert
                s
            \Vert_{Y(T)}
        ).
    \end{gather*}
\end{proposition}
\begin{proof}
    By the assumptions, we see from \eqref{eq_formula_dw_2} that
    \begin{align*}
        & \vert
            \partial_y w_2(t)
        \vert\\
        & \leq \Vert
            f
        \Vert_{Lip(\Real \times (0, T))}
        \int_0^t
            \vert
                \partial_y G_{t - \tau} (y + s(t) - s(\tau))
            \vert
            (
                \vert
                      s(\tau)
                    - s(t)
                \vert
                + \vert
                      \tau
                    - t
                \vert
            )
        d\tau\\
        & + \Vert
            f
        \Vert_{BUC(\Real \times (0, T))}
        \left \vert
            \int_0^t
                \partial_y G_{t - \tau} (y + s(t) - s(\tau))
            d\tau.
        \right \vert
    \end{align*}
    
    The integrand in the first term is bounded by $C(\Vert s \Vert_{Y(T)} + 1)$.
    The integral in the second term is bounded \red{in a similar way to that for} $w_3$; see \eqref{eq_int_dG_tau_dtau_I1} and \eqref{eq_int_dG_tau_dtau_I2} in the proof of Proposition \ref{prop_estimates_for_w3}.
    Therefore, we have
    \begin{align*}
        \sup_{0 < t < T}
            \Vert
                \partial_y w_2(t)
            \Vert_\infty
        & \leq C T
        \Vert
            f
        \Vert_{BUC((0, T); Lip(\Real))}
        (
            \Vert
                s
            \Vert_{Y(T)}
            + 1
        )\\
        & + C \Vert
            f
        \Vert_{BUC(\Real \times (0, T))}
        (
            1
            + T^\frac{1}{2}
            \Vert
                s
            \Vert_{Y(T)}
        ).
    \end{align*}
\end{proof}

\begin{proposition} \label{prop_estimates_for_w2_for_smooth_f_diff}
    Let $s_1, s_2 \in Y(T)$, and $f \in BUC((0, T); BUC^1(\Real))$.
    Assume that $T < 1$.
    Then
    \begin{align*}
        & \sup_{0 < t < T}
            \vert
                \gamma_0 \partial_y
                (
                    w_2(s_2; t)
                    - w_2(s_1; t)
                )
            \vert\\
        & \leq C T^\frac{1}{2}
        \Vert
            f
        \Vert_{BUC((0, T); BUC^1(\Real))}
        (
            1
            + \Vert
                s_1
            \Vert_{Y(T)}
        )
        \Vert
            \sigma
        \Vert_{Y(T)}.
    \end{align*}
\end{proposition}
\begin{proof}
    Using the formula \eqref{eq_formula_gamma_0_dw2s2_w2s1}, we deduce that
    \begin{align*}
        & \vert
            \gamma_0 \partial_y
            (
                  w_2(s_2; t)
                - w_2(s_1; t)
            )
        \vert\\
        & \leq C
        \Vert
            f
        \Vert_{BUC(\Real \times (0, T))}
        \Vert
            \sigma
        \Vert_{Y(T)}
        \int_0^t
            (
                t
                - \tau
            )^{- \frac{1}{2}}
        d\tau\\
        & + C \Vert
            \partial_x f
        \Vert_{BUC(\Real \times (0, T))}
        \int_0^t
            (
                t
                - \tau
            )^{- \frac{1}{2}}
            \tau
        d\tau
        \Vert
            s_1
        \Vert_{Y(T)}
        \Vert
            \sigma
        \Vert_{Y(T)}\\
        & \leq C
        (
            T^\frac{1}{2}
            \Vert
                f
            \Vert_{BUC(\Real \times (0, T))}
            + C T^{\frac{3}{2}}
            \Vert
                \partial_x f
            \Vert_{BUC(\Real \times (0, T))}
            \Vert
                s_1
            \Vert_{Y(T)}
        )
        \Vert
            \sigma
        \Vert_{Y(T)}.
    \end{align*}
\end{proof}

\subsection{Estimates for $w_4$} \label{subsec_w4}
We finally estimate the forcing term $w_4$.
%
%\begin{align*}
%    w_4(t)
%    = \int_0^{t}
%        \int_0^\infty
%            G_{t - \tau} (y - z + s(t) - s(\tau))
%            \partial_x f(z + s(\tau), \tau)
%        dz
%    d\tau
%\end{align*}
%
\begin{proposition} \label{prop_estimates_for_w4}
    Let $q \in [1, \infty]$, $r \in (2, \infty]$, and $\eta \in [0, 1 - \frac{1}{r})$.
    Let $s \in C[0, T]$.
    Let $f \in L^1_{loc}(0, T; L^1_{loc}(\Real))$ satisfy
    \begin{align*}
        t^{\eta} \partial_x f
        \in L^r(0, T; L^q(\Real)).
    \end{align*}
    Then
    \begin{align*}
        \sup_{0 < t < T}
            \Vert
                w_4(t)
            \Vert_{L^q(\Omega)}
        + \sup_{0 < t < T}
            t^{\frac{1}{2}}
            \Vert
                \partial_y w_4(t)
            \Vert_{L^q(\Omega)}
        \leq C
        \Vert
            t^{\eta} \partial_x f
        \Vert_{L^r(0, T; L^q(\Real))}.
    \end{align*}
    The constant $C > 0$ is $O(T^{1 - \frac{1}{r} - {\eta}})$.
\end{proposition}
\begin{proof}
    We find from the Young inequality and Proposition \ref{prop_beta_type} that
    \begin{align*}
        \Vert
            \partial_y^m w_4(t)
        \Vert_{L^q(\Omega)}
        & \leq C
        \int_0^t
            (t - \tau)^{- \frac{m}{2}}
            \Vert
                \partial_x f(\tau)
            \Vert_{L^q(\Real)}
        d\tau\\
        & \leq C t^{1 - \frac{1}{r} - \frac{m}{2} - {\eta}}
        \Vert
            t^{\eta} \partial_x f
        \Vert_{L^r(0, T; L^q(\Real))}
    \end{align*}
    for $m = 0,1$.
\end{proof}

\begin{proposition} \label{prop_estimates_for_diff_w4}
    Assume the same assumptions as \red{in} Proposition \ref{prop_estimates_for_w4} for $f, T, q, r, \eta$.
    Let $s_1, s_2 \in X^\frac{1}{2}(T)$.
    Assume further
    \begin{align*}
        t^{\frac{1}{2} + \eta} \partial_x^2 f
        \in L^r(0, T; L^q(\Real)).
    \end{align*}
    Then
    \begin{align*}
        & \sup_{0 < t < T}
            t^{\frac{1}{2}}
            \vert
                \gamma_0 \partial_y
                (
                      w_4(s_2; t)
                    - w_4(s_1; t)
                )
            \vert\\
        & \leq C
        \Vert
            \sigma
        \Vert_{X^\frac{1}{2}(T)}
        (
            \Vert
                t^{\eta} \partial_x f
            \Vert_{L^r(0, T; L^q(\Real))}
            + \Vert
                t^{\frac{1}{2} + \eta} \partial_x^2 f
            \Vert_{L^r(0, T; L^q(\Real))}
        ).
    \end{align*}
    The constant $C > 0$ is $O(T^{1 - \frac{1}{r} - \eta + \frac{1}{2q}})$.
\end{proposition}
\begin{proof}
    We first derive the formula
    \begin{align} \label{eq_formula_w4s2_w4s1}
        \begin{split}
            & w_4(s_2; t)
            - w_4(s_1; t)\\
            & = \int_0^t
                \int_0^\infty
                    (
                          \sigma(t)
                        - \sigma(\tau)
                    )
                    \int_0^1
                        \partial_y G_{t - \tau}
                        (
                            y
                            - z
                            + s_\mu(t)
                            - s_\mu(\tau)
                        )
                    d\mu
                    \partial_z f(z + s_2(\tau), \tau)
                dz
            d\tau\\
            & + \int_0^{t}
                \int_0^\infty
                    G_{t - \tau} (y - z + s_1(t) - s_1(\tau))
                    \sigma(\tau)
                    \int_0^1
                          \partial_z^2 f(z + s_\mu(\tau), \tau)
                    d\mu
                dz
            d\tau.
        \end{split}
    \end{align}
    We deduce from Proposition \ref{prop_beta_type} that
    \begin{align*}
        & \vert
            \gamma_0 \partial_y
            (
                  w_4(s_2; t)
                - w_4(s_1; t)
            )
        \vert\\
        & \leq C
        \int_0^t
            \vert
                  \sigma(t)
                - \sigma(\tau)
            \vert
            (t - \tau)^{- \frac{1}{2} - \frac{1}{2q}}
            \Vert
                \partial_x f(\tau)
            \Vert_{L^q(\Real)}
        d\tau\\
        & + C
        \int_0^t
            \vert
                \sigma(\tau)
            \vert
            (
                t
                - \tau
            )^{- \frac{1}{2q}}
            \Vert
                \partial_x^2 f(\tau)
            \Vert_{L^q(\Real)}
        d\tau\\
        & \leq C t^{1 - \frac{1}{r} + \frac{1}{2q} - \eta}
        \Vert
            \sigma
        \Vert_{X^\frac{1}{2}(T)}\\
        & \quad\quad\quad\quad \times
        \Bigl(
            \Vert
                t^{\eta} \partial_x f
            \Vert_{L^r(0, T; L^q(\Real))}
            + \Vert
                t^{\frac{1}{2} + \eta} \partial_x^2 f
            \Vert_{L^r(0, T; L^q(\Real))}
        \Bigr).
    \end{align*}
\end{proof}
We consider the case \red{where} $s \in Y(T)$.
The following proposition is a direct consequence of Propositions \ref{prop_estimates_for_w4} and \ref{prop_estimates_for_diff_w4} \red{with} $f$ \red{replaced} by $\partial_x f$.
\begin{proposition} \label{prop_estimates_for_w4_for_smooth_f}
    Let $s \in Y(T)$ and $f \in BUC((0, T); BUC^2(\Real))$.
    Assume that $T < 1$.
    Then
    \begin{gather*}
        \sup_{0 < t < T}
            \Vert
                \partial_y w_4(t)
            \Vert_{\infty}
        \leq C
        \Vert
            f
        \Vert_{BUC((0, T); BUC^2(\Real))}
        (
            1
            + T^\frac{1}{2}
            \Vert
                s
            \Vert_{Y(T)}
        ).
    \end{gather*}
\end{proposition}
Similarly, we have
\begin{proposition} \label{prop_estimates_for_w4_for_smooth_f_diff}
    Let $s_1, s_2 \in Y(T)$ and $f \in BUC((0, T); BUC^2(\Real))$.
    Assume that $T < 1$.
    Then
    \begin{align*}
        & \sup_{0 < t < T}
            \vert
                \gamma_0 \partial_y
                (
                      w_4(s_2; t)
                    - w_4(s_1; t)
                )
            \vert\\
        & \leq C T^\frac{1}{2}
        \Vert
            f
        \Vert_{BUC((0, T); BUC^2(\Real))}
        (
            1
            + \sum_{j=1,2}
                \Vert
                    s_j
                \Vert_{Y(T)}
        )
        \Vert
            \sigma
        \Vert_{Y(T)}.
    \end{align*}
\end{proposition}
\begin{proof}
    \red{Replacing} $f$ \red{by} $\partial_x f$ in Proposition \ref{prop_estimates_for_w2_for_smooth_f_diff} and $\partial_x f$ \red{by} $\partial_x^2 f$ in the formula \eqref{eq_formula_w4s2_w4s1} and using the Young inequality, we have
    \begin{align*}
        & \vert
            \gamma_0 \partial_y
            (
                  w_4(s_2; t)
                - w_4(s_1; t)
            )
        \vert\\
        & \leq C T^\frac{1}{2}
        \Vert
            f
        \Vert_{BUC((0, T); BUC^2(\Real))}
        (
            1
            + \sum_{j=1,2}
                \Vert
                    s_j
                \Vert_{Y(T)}
        )
        \Vert
            \sigma
        \Vert_{Y(T)}\\
        & + C T^\frac{1}{2}
        \Vert
            \sigma
        \Vert_{Y(T)}
        \Vert
            f
        \Vert_{BUC((0, T); BUC^2(\Real))}.
    \end{align*}
\end{proof}

\subsection{Construction of the solution}

We construct the solution $s$ in $X^\frac{1}{2}(T)$ and $Y(T)$ by Banach's fixed point theorem.

It \red{is} worth \red{mentioning} that Lemma \ref{lem_local_existence_of_w} below \red{is} valid without Assumption (C).
A mild solution \red{obtained in} Lemma \ref{lem_local_existence_of_w} is a candidate \red{for a} solution to \eqref{eq_w}, which might not satisfy $\gamma_0 w(t) = \alpha$.
Assumption (C) ensures the monotonicity of $s$ at least \red{for a} short time, \red{and} the monotonicity \red{leads} to $\gamma_0 w(t) = \alpha$, which is proved in the next section.

We denote by
\begin{align*}
    N(s, t)
    := - \frac{\gamma_0 \partial_y w(s)}{\alpha}
    - \frac{f(s(t), t)}{\alpha}.
\end{align*}
Summing up \red{the} estimates in \red{Sections} \ref{subsec_w1_w3}\red{--}\ref{subsec_w4}, we obtain
\begin{corollary}\label{cor_estimates_for_w}
    Let $s, s_1, s_2 \in X^\frac{1}{2}(T)$ and $w_0 \in BUC(\Omega)$ satisfying $\gamma_0 w_0 = \alpha$.
    %Let $\varepsilon \in (0, \frac{1}{2})$.
    Let $w_0$ and $f$ satisfy the same assumptions \red{as in} Theorem \ref{thm_main1}.
    Then
    \begin{align} \label{eq_estimtate_N}
        \begin{split}
            & \sup_{0 < t < T}
                t^\frac{1}{2}
                \vert
                    N(s(t), t)
                \vert\\
            %& \leq
            %\left(
            %    C_{w_0, T}
            %    + C_1 \alpha
            %        \Vert
            %        s
            %    \Vert_{X^\theta(T)}^2
            %\right)\\
            %& + C_2 T^\frac{1}{2}
            %\Vert
            %    t^\frac{1}{2} \partial_x f
            %\Vert_{L^\infty(0, T; L^\infty(\Real))}\\
            %& + C_3 C_f \left(
            %    1
            %    + \Vert
            %        s
            %    \Vert_{X^\frac{1}{2}(T)}
            %\right)\\
            & \leq \frac{C_0}{\alpha}
            (
                C_{w_0, T}
                + C_f
                + \Vert
                    t^\frac{1}{2} f
                \Vert_{BUC(\Real \times (0, T))}
                + T^{\varepsilon}
                \Vert
                    t^{1 - \varepsilon} \partial_x f
                \Vert_{BUC(\Real \times (0, T))}
            )\\
            & + \frac{
                C_1
                (
                    C_f
                    + 
                    \Vert
                        t^\frac{1}{2} f
                    \Vert_{BUC(\Real \times (0, T))}
                )
            }{\alpha}
            \Vert
                s
            \Vert_{X^\frac{1}{2}(T)}
            + C_2
            \Vert
                s
            \Vert_{X^\frac{1}{2}(T)}^2
        \end{split}
    \end{align}
    and
    \begin{align} \label{eq_estimtate_diff_N}
        \begin{split}
            & \sup_{0 < t < T}
                t^{\frac{1}{2}}
                \vert
                      N(s_2(t), t)
                    - N(s_1(t), t)
                \vert\\
            & \leq
            \left(
                \frac{C_{w_0, T}^\prime}{\alpha}
                + C_1^\prime
                \sum_{j=1,2}
                    \Vert
                        s_j
                    \Vert_{X^\frac{1}{2}(T)}
            \right)
            \Vert
                \sigma
            \Vert_{X^\frac{1}{2}(T)}\\
            & + \frac{C_2^\prime}{\alpha}
            \Biggl(
                \tnorm{f}_{1, T}
                + (
                    T^\frac{1}{2}
                    + \sum_{j=1,2}
                        \Vert
                            s_j
                        \Vert_{X^\frac{1}{2}(T)}
                    \tnorm{f}_{2, T}
                )
            \Biggr)
            \Vert
                \sigma
            \Vert_{X^\frac{1}{2}(T)},
        \end{split}
    \end{align}
    where
    \begin{align*}
        & \tnorm{f}_{1, T}
        =
        \Vert
            t^\frac{1}{2} f(t)
        \Vert_{BUC(\Real \times (0, T))}
        + \tnorm{f}_{2, T},\\
        & \tnorm{f}_{2, T}
        = \sum_{m=1,2}
            \Vert
                t^{\frac{m+1}{2} - \varepsilon} \partial_y^m f(t)
            \Vert_{BUC(\Real \times (0, T))},
    \end{align*}
    and the constants $C_j, C_j^\prime > 0$ are independent of $T$.
\end{corollary}
\begin{corollary}\label{cor_estimates_for_w_when_w0_f_are_smooth}
    Let $T < 1$.
    Let $s, s_1, s_2 \in Y(T)$.
    Let $w_0$ and $f$ satisfy the same assumptions \red{as in} Theorem \ref{thm_main2}.
    Then
    \begin{align*}
        \sup_{0 < t < T}
            \vert
                \gamma_0\partial_y w(s; t)
            \vert
        %& \leq C
        %\left(
        %    \Vert
        %        \partial_y w_0
        %    \Vert_{\infty}
        %    + \alpha  T^\frac{1}{2}
        %    \Vert
        %        s
        %    \Vert_{Y(T)}^2
        %\right)\\
        %& + C \Vert
        %    f
        %\Vert_{BUC^1(Q_T)}
        %(
        %    1
        %    + T^\frac{1}{2}
        %    \Vert
        %        s
        %    \Vert_{Y(T)}
        %)\\
        %& + C
        %\Vert
        %    f
        %\Vert_{BUC^2(Q_T)}
        %(
        %    1
        %    + T^\frac{1}{2}
        %    \Vert
        %        s
        %    \Vert_{Y(T)}
        %)\\
        & \leq C_0
        \left(
            \Vert
                \partial_y w_0
            \Vert_{\infty}
            + \Vert
                f
            \Vert_{BUC((0, T); BUC^2(\Real))}
        \right)\\
        & + C_1 T^\frac{1}{2} \Vert
            f
        \Vert_{BUC((0, T); BUC^2(\Real))}
        \Vert
            s
        \Vert_{Y(T)}
        + C_2 \alpha T^\frac{1}{2}
            \Vert
                s
            \Vert_{Y(T)}^2
    \end{align*}
    and
    \begin{align*}
        & \sup_{0 < t < T}
            \vert
                \gamma_0\partial_y
                (
                    w(s_2; t)
                    - w(s_1; t)
                )
            \vert\\
        %& \leq C T^\frac{1}{2}
        %\left(
        %    \Vert
        %        \partial_y w_0
        %    \Vert_{\infty}
        %    + \alpha
        %    \sum_{j=1,2}
        %        \Vert
        %            s_j
        %        \Vert_{Y(T)}
        %\right)
        %\Vert
        %    \sigma
        %\Vert_{Y(T)}\\
        %& + C T^\frac{1}{2}
        %\Vert
        %    f
        %\Vert_{BUC^1(Q_T)}
        %(
        %    1
        %    + \sum_{j=1,2}
        %        \Vert
        %            s_j
        %        \Vert_{Y(T)}
        %)
        %\Vert
        %    \sigma
        %\Vert_{Y(T)}\\
        %& + C T^\frac{1}{2}
        %(
        %    \Vert
        %        f
        %    \Vert_{BUC^2(Q_T)}
        %    + \Vert
        %        f
        %    \Vert_{L^1(0, T; BUC^3(\Omega))}
        %)
        %(
        %    1
        %    + \sum_{j=1,2}
        %        \Vert
        %            s_j
        %        \Vert_{Y(T)}
        %)
        %\Vert
        %    \sigma
        %\Vert_{Y(T)}\\
        & \leq C T^\frac{1}{2}
        \left(
            \Vert
                \partial_y w_0
            \Vert_{\infty}
            + \alpha
            \sum_{j=1,2}
                \Vert
                    s_j
                \Vert_{Y(T)}
        \right)
        \Vert
            \sigma
        \Vert_{Y(T)}\\
        & + C T^\frac{1}{2}
        \Vert
            f
        \Vert_{BUC((0, T); BUC^2(\Real))}
        (
            1
            + \sum_{j=1,2}
                \Vert
                    s_j
                \Vert_{Y(T)}
        )
        \Vert
            \sigma
        \Vert_{Y(T)}.
    \end{align*}
\end{corollary}
Corollary \ref{cor_estimates_for_w} implies that the right-hand side of
\begin{align*}
    \dot{s}
    %= - \frac{\gamma_0 \partial_y w(s)}{\alpha}
    %- \frac{f(s(t), t)}{\alpha}
    = N(s, t)
\end{align*}
is locally Lipschitz continuous in $X^\frac{1}{2}(T)$.
If $C_f \leq \frac{\alpha}{2 C_1}$ and if $T > 0$ and the first term on the right-hand side of \eqref{eq_estimtate_N} is sufficiently small, we find that the solution mapping $s \mapsto \int_0^t N(s(\tau), \tau) d\tau$ in $X^\frac{1}{2}(T)$ maps the ball of radius $r$ into itself.
The radius is bounded by
\begin{align*}
    r \leq \frac{
        4 C_0 (
            C_{w_0, T}
            + C_f
            + T^\varepsilon
                \Vert
                    t^{1 - \varepsilon} \partial_x f
                \Vert_{BUC(\Real \times (0, T))}
        )
    }{\alpha}.
\end{align*}
The estimate \eqref{eq_estimtate_diff_N} implies that $s \mapsto \int_0^t N(s(\tau), \tau) d\tau$ is a contraction mapping if $r$, $T$, and $f$ are sufficiently small.
Therefore, we obtain a unique solution $s$ to
\begin{align} \label{eq_s_prime_equal_N}
    \dot{s}
    = N(s, \cdot)
\end{align}
by Banach's fixed point theorem \red{on} some short interval $[0, T_\ast]$.
Once $s$ is given, $w$ is automatically determined by \eqref{eq_mild_solution_w}.
We obtain
\begin{lemma} \label{lem_local_existence_of_w}
    Assume the same hypotheses as in Corollary \ref{cor_estimates_for_w}.
    Assume further
    \begin{align*}
        w_0 \in L^1(\Omega), \quad
        t^{1 - \varepsilon} \partial_x f \in L^\infty(0, T; L^1(\Omega))
    \end{align*}
    for some $\varepsilon \in (0, \frac{1}{2})$.
    Then there exist $T_\ast \in (0, T)$ and $\delta > 0$ such that, if
    \begin{align*}
        C_{w_0, T}
        + C_f
        + \tnorm{f}_{1, T}
        < \delta,
    \end{align*}
    there exists a unique mild solution $(s, w) \in X^\frac{1}{2}(T) \times C([0, T_\ast); L^1(\Omega) \cap BUC(\Omega))$ to \eqref{eq_w} such that
    \begin{align} \label{eq_estimate_for_s_lem_local_existence_of_w}
        & \alpha \Vert
            s
        \Vert_{X^\frac{1}{2}(T_\ast)}
        \leq C \delta
    \end{align}
    and
    \begin{align} \label{eq_estimate_for_w_lem_local_existence_of_w}
        &\sup_{0 < t < T_\ast}
            \Vert
                w(t)
            \Vert_{L^q(\Omega)}
        + \sup_{0 < t < T_\ast}
            t^{\frac{1}{2}}
            \Vert
                \partial_y w(t)
            \Vert_{L^q(\Omega)}%\\
        %& \leq C \Vert
        %    w_0
        %\Vert_{L^{q}(\Omega)}\\
        %& + \alpha C
        %\left(
        %    1
        %    +
        %    \Vert
        %        s
        %    \Vert_{X^\frac{1}{2}(T)}
        %    +
        %    \Vert
        %        s
        %    \Vert_{X^\frac{1}{2}(T)}^2
        %\right)\\
        %& + C T^\frac{1}{2} \Vert
        %    t^\frac{1}{2} f
        %\Vert_{L^\infty(0, T; L^q(\Omega))}\\
        %& + C \Vert
        %    t^{\frac{1}{2}} f
        %\Vert_{L^\infty(0, T; L^{\infty}(\Real))}
        %+ C C_f
        %\left(
        %    1
        %    + \Vert
        %        s
        %    \Vert_{X^\frac{1}{2}(T)}
        %\right)\\
        < \infty,
    \end{align}
    for all $q \in [1, \infty]$.
\end{lemma}
\begin{remark}
    Lemma \ref{lem_local_existence_of_w} does not ensure that $\gamma_0 w(t) \equiv \alpha$.
\end{remark}
\begin{remark}
    In our previous work, we have presented \red{a} bounded self-similar solution with zero initial data; see Theorem 7 in \cite{FurukawaGigaKajiwara2025}.
    This observation does not contradict the previous result.
    Indeed, since we assume that $\gamma_0 w_0 = \alpha$, such self-similar solutions have already been excluded.
\end{remark}
\begin{remark} %\label{rmk_global_existence_of_w}
    We observe \eqref{eq_estimtate_N} again.
    When
    \begin{align} \label{eq_simple_f}
        f(x, t)
        = \frac{f_0}{\sqrt{t}}
    \end{align}
    for some positive constant $f_0$, the derivative $\partial_x f$ and $C_f$ vanish.
    Assume that
    \begin{align} \label{eq_closeness_of_trivial_sol}
        \left \Vert
            \frac{w_0}{\alpha}
            - 1
        \right \Vert_{\infty}
        + \frac{f_0}{\alpha}
        \leq \delta
    \end{align}
    for some sufficiently small $\delta > 0$.
    Repeating the same argument \red{as in the proof of} Lemma \ref{lem_local_existence_of_w}, we find from \eqref{eq_estimtate_N} and \eqref{eq_estimtate_diff_N} that the mapping $s \mapsto \int_0^t N(s(\tau), \tau) d\tau$ is a contraction self-mapping in $X^\frac{1}{2}(T)$ for all $T > 0$.
    In this case, we obtain \red{a} global solution $(w, s)$ satisfying the same estimate as in Lemma \ref{lem_local_existence_of_w}, and the corresponding $w$ may fail to satisfy the boundary condition.
    %Assume the same assumptions for $w_0$ as Corollary \ref{cor_estimates_for_w}.
    %Assume that $f$ is of the form \eqref{eq_simple_f} for some constant $f_0 \geq 0$.
    Under this assumption, $v_0(y) := \int_0^y w_0(z)\,dz$ and $v(y, t) := \int_0^y w(z, t)\,dz$ are necessarily unbounded.
    As an example of an unbounded solution with $f \equiv 0$, we have
    \begin{align*}
        v(y, t)
        = \alpha y, \quad
        s(t)
        \equiv 0, \quad
        \text{for all $t\geq 0$}.
    \end{align*}
    The condition \eqref{eq_closeness_of_trivial_sol} can be rewritten as the existence of small global solutions around the trivial solution.
    However, we do not consider the boundary values or boundary behavior of such solutions in this paper.
    This is because, in the main results of this paper, the boundary is required to be moving in order for the solution to satisfy the overdetermined boundary condition.
    %Consequently, solutions in a neighborhood of stationary solutions to the overdetermined boundary value problem are excluded.
\end{remark}
Using Corollary \ref{cor_estimates_for_w_when_w0_f_are_smooth} and the same argument as \red{in} Lemma \ref{lem_local_existence_of_w}, we obtain
\begin{lemma} \label{lem_local_existence_of_w_when_w0_f_are_smooth}
    Assume the same assumptions \red{as in} Theorem \ref{thm_main2} for $w_0$ and $f$.
    Assume further
    \begin{align*}
        w_0 \in L^1(\Omega), \quad
        f \in BUC((0, T); L^1(\Omega)).
    \end{align*}
    Then there exist $T_\ast \in (0, T)$ and a unique solution $(s, w) \in Y(T_\ast) \times C([0, T_\ast); L^1(\Omega) \cap BUC^1(\Omega))$ to \eqref{eq_w} such that
    \begin{align} %\label{eq_estimate_for_s_lem_local_existence_of_w}
        & \alpha \Vert
            s
        \Vert_{Y(T_\ast)}
        +
        \sup_{0 < t < T_\ast}
            \Vert
                w(t)
            \Vert_{BUC^1(\Omega)}
        < \infty.
    \end{align}
\end{lemma}
\begin{remark}
    For the solution obtained in Lemma \ref{lem_local_existence_of_w_when_w0_f_are_smooth}, the quantity $\gamma_0 \partial_y w + f(s, \cdot)$ is continuous over $[0, T]$.
    Therefore\red{,} $s$ is monotone at least on a short interval provided that $\gamma_0 \partial_y w_0 + f(0, 0) \neq 0$.
\end{remark}

\section{Boundary condition $\gamma_0 \partial_y v = \alpha$} \label{sec_boundary_condition}
It remains to show that $w$ satisfies the boundary condition $\gamma_0 w(t) = \alpha$ for all $t \geq 0$.
This corresponds to the contact angle condition for $v$.
If $s \equiv 0$, the problem is overdetermined.
In this case, the mild solution constructed in Section~ \ref{sec_mild_solution} cannot attain the prescribed boundary value $\alpha$.
Otherwise, provided that no strong oscillation occurs near $t=0$, the mild solution satisfies the boundary condition.
\begin{lemma} \label{lem_Neumann_boundary_condition}
    Assume (A) and (B).
    Let $v_0 \in BUC^1(\Omega)$ be a positive function satisfying $\partial_y v_0 \in L^1(\Omega)$, $\gamma_0 v_0 = 0$, and $\gamma_0 \partial_y v_0 = \alpha$.
    Let
    \begin{align*}
        (w, s)
        \in C([0, T); BUC(\Omega))
        \cap C((0, T); BUC^1(\Omega))
        \times X^\frac{1}{2}(T)
    \end{align*}
    be the mild solution.
    Assume that $s$ is strictly monotone on $[0, T_\ast]$ for some $T_\ast < T$.% and $\dot{s}(t) \neq 0$ for all $t \in (0, T_\ast]$.
    Then
    \begin{align*}
        \gamma_0 w
        \equiv \alpha
    \end{align*}
    on $[0, T_\ast]$.
\end{lemma}
\begin{proof}
    %Let us denote the approximation of $s$ by $s(t) = \frac{1}{\varepsilon} \int_t^{t+\varepsilon} s(\tau) d\tau \in Y(T_\varepsilon)$ for $T_\varepsilon = T - \varepsilon$.
    Let $\varepsilon \in (0, 1)$ be a small number.
    We consider the approximating solution $(W_\varepsilon(t), S_\varepsilon(t))$ defined by
    \begin{align*}
        \begin{split}
            W_\varepsilon(y, t)
            & := w(\varepsilon + y, t)\\
            & = \int_0^\infty
                G_t(
                    \varepsilon
                    + y
                    - z
                    + s(t)
                )
                w_0(z)
            dz\\
            & + \int_0^{t}
                G_{t - \tau}
                (
                    \varepsilon
                    + y
                    + s(t)
                    - s(\tau)
                )
                f(s(\tau), \tau)
            d\tau\\
            & - \int_0^{t}
                \alpha \partial_y G_{t - \tau}
                (
                    \varepsilon
                    + y
                    + s(t)
                    - s(\tau)
                )
            d\tau\\
            & + \int_0^{t}
                \int_{0}^\infty
                    G_{t - \tau}
                    (
                        \varepsilon
                        + y
                        - z
                        + s(t)
                        - s(\tau)
                    )
                    \partial_z f(z + s(\tau), \tau)
                dz
            d\tau%\\
            %& =: N_\varepsilon(w_0, f; s),
        \end{split}
    \end{align*}
    and
    \begin{align*}
        S_\varepsilon(t)
        := - \frac{1}{\alpha}
        \int_0^t
            \partial_y W_\varepsilon (0, \tau)
            + f (s(\tau), \tau)
        d\tau, \quad 
        a.a. \; t \in (0, T_\varepsilon).
    \end{align*}
    In view of the estimates obtained in Sections \ref{subsec_w1_w3}--\ref{subsec_w4}, $(W_\varepsilon,S_\varepsilon)$ belongs to the same solution space as $(w,s)$, uniformly in $\varepsilon > 0$.
    Moreover, by construction, $(W_\varepsilon,S_\varepsilon)$ is smooth, since the approximating integral kernels are smooth for $y \geq 0$ and $t \geq 0$.
    The function $W_\varepsilon$ solves
    \begin{align*}
        \partial_t W_\varepsilon
        - \partial_y^2 W_\varepsilon
        &= \dot{s} \partial_y W_\varepsilon
        + \partial_y f(\cdot + \varepsilon + s(t), t)\\
        &=: \dot{s} \partial_y W_\varepsilon
        + \partial_y f_\varepsilon(\cdot + s(t), t),
    \end{align*}
    for $y > 0$ and $t \in (0, T_\varepsilon)$ in \red{the} classical sense.

    It is clear from the definition of $S_\varepsilon$ that
    \begin{align} \label{eq_def_J1_J2}
        \begin{split}
            0
            &= 
            - \left(
                  f(s(\tau), \tau)
                - f_\varepsilon(s(\tau), \tau)
            \right)
            G_{t - \tau} (
                \varepsilon
                + y
                + s(t)
                - s(\tau)
            )\\
            & - \left(
                \alpha \dot{S_\varepsilon}(\tau)
                + f_\varepsilon(s(\tau), \tau)
                + \partial_y W_\varepsilon(0, \tau)
            \right)
            G_{t - \tau} (
                \varepsilon
                + y
                + s(t)
                - s(\tau)
            )\\
            & =: J_1 + J_2
        \end{split}
    \end{align}
    then by the fundamental theorem \red{of} calculus
    \begin{align*}
        J_2
        & = - \int_0^\infty
            \partial_z
            \left[
                \left(
                    \alpha \dot{S_\varepsilon}(\tau)
                    + f_\varepsilon(z + s(\tau), \tau)
                    + \partial_y W_\varepsilon(z, \tau)
                \right)
            \right.\\
        & \quad\quad\quad\quad\quad\quad\quad\quad
        \times
            \left.
                G_{t - \tau} (
                    \varepsilon
                    + y
                    - z
                    + s(t)
                    - s(\tau)
                )
            \right]
        dz
    \end{align*}
    Let $\delta \in [\varepsilon, T)$.
    Integrating over $(\delta, t - \delta)$, we deduce that
    \begin{align} \label{eq_integration_wrt_ytau_boundary_ode}
        \begin{split}
                \int_\delta^{t - \delta}
                    J_2
                d\tau
            & = - \int_{\delta}^{t - \delta}
                (
                    \alpha
                    - W_\varepsilon(0, \tau)
                )
                \dot{S_\varepsilon}(\tau)
                G_{t - \tau}
                (
                    \varepsilon
                    + y
                    + s(t)
                    - s(\tau)
                )
            d\tau\\
            & - \int_{\delta}^{t - \delta}
                W_\varepsilon(0, \tau)
                (
                    \dot{S_\varepsilon}(\tau)
                    - \dot{s}(\tau)
                )
                G_{t - \tau}
                (
                    \varepsilon
                    + y
                    + s(t)
                    - s(\tau)
                )
            d\tau\\
            & + \int_{\delta}^{t - \delta}
                \int_{0}^\infty
                    \partial_z
                    \bigr [
                        \{
                            W_\varepsilon(z, \tau) \dot{s}(\tau)
                            + f_\varepsilon(z + s(\tau), \tau)\\
            & \quad\quad\quad\quad\quad\quad
                            + \partial_z W_\varepsilon(z, \tau)
                        \}
                        G_{t - \tau}
                        (
                            \varepsilon
                            + y
                            - z
                            + s(t)
                            - s(\tau)
                        )
                    \bigl ]
                dz
            d\tau.
        \end{split}
    \end{align}
    We expand the integrand in the last term by using the equation that $W_\varepsilon$ \red{satisfies}.
    Then the last term becomes
    \begin{align} \label{eq_int_int_dtW_G}
        \begin{split}
            & \int_{\delta}^{t - \delta}
                \int_{0}^\infty
                    \partial_\tau W_\varepsilon(z, \tau)
                    G_{t - \tau}
                    (
                        \varepsilon
                        + y
                        - z
                        + s(t)
                        - s(\tau)
                    )
                dz
            d\tau\\
            & - \int_{\delta}^{t - \delta}
                \int_{0}^\infty
                    \left(
                        W_\varepsilon(z, \tau) \dot{s}(\tau)
                        + f_\varepsilon(s(\tau) + z, \tau)
                        + \partial_z W_\varepsilon(z, \tau)
                    \right)\\
            & \quad\quad\quad\quad\quad \times
                    \partial_y G_{t - \tau}
                    (
                        \varepsilon
                        + y
                        - z
                        + s(t)
                        - s(\tau)
                    )
                dz
            d\tau.
        \end{split}
    \end{align}
    Here we derive the two formulas
    \begin{align} \label{eq_sec3_formula_for_integrand_1}
        \begin{split}
            & \partial_\tau W_\varepsilon(z, \tau)
            G_{t - \tau}
            (
                \varepsilon
                + y
                - z
                + s(t)
                - s(\tau)
            )\\
            & = \partial_\tau
            \left[
                W_\varepsilon(z, \tau)
                G_{t - \tau}
                (
                    \varepsilon
                    + y
                    - z
                    + s(t)
                    - s(\tau)
                )
            \right]\\
            & + W_\varepsilon(z, \tau)
            \partial_t G_{t - \tau}
            (
                \varepsilon
                + y
                - z
                + s(t)
                - s(\tau)
            )\\
            & + 
                W_\varepsilon(z, \tau) \dot{s}(\tau) \partial_y G_{t - \tau}
                (
                    \varepsilon
                    + y
                    - z
                    + s(t)
                    - s(\tau)
                )
            .
        \end{split}
    \end{align}
    and, by integration by parts,
    \begin{align} \label{eq_sec3_formula_for_integrand_2}
        \begin{split}
            & \int_{\delta}^{t - \delta}
                \int_0^\infty
                    [
                        f_\varepsilon(s(\tau) + z, \tau)
                        + \partial_z W_\varepsilon(z, \tau)
                    ]\\
            & \quad\quad\quad\quad\quad\quad
                \times
                    \partial_y G_{t - \tau}
                    (
                        \varepsilon
                        + y
                        - z
                        + s(t)
                        - s(\tau)
                    )
                dz
            d\tau\\
            & = \int_{\delta}^{t - \delta}
                f_\varepsilon(s(\tau), \tau)
                G_{t - \tau}
                (
                    \varepsilon
                    + y
                    + s(t)
                    - s(\tau)
                )
            d\tau\\
            & + \int_{\delta}^{t - \delta}
                \int_0^\infty
                    \partial_z f_\varepsilon(z + s(\tau), \tau)
                    G_{t - \tau}
                    (
                        \varepsilon
                        + y
                        - z
                        + s(t)
                        - s(\tau)
                    )
                dz
            d\tau\\
            & - \int_{\delta}^{t - \delta}
                W_\varepsilon(0, \tau)
                \partial_y G_{t - \tau}
                (
                    \varepsilon
                    + y
                    %- z
                    + s(t)
                    - s(\tau)
                )
            d\tau\\
            & + \int_{\delta}^{t - \delta}
                \int_0^\infty
                    W_\varepsilon(z, \tau)
                    \partial_y^2 G_{t - \tau}
                    (
                        \varepsilon
                        + y
                        - z
                        + s(t)
                        - s(\tau)
                    )
                dz
            d\tau.
        \end{split}
    \end{align}
    Using the equation satisfied by $G_{t-\tau}(\varepsilon + y - z + s(t) - s(\tau))$, the space--time integrals of the right-hand side of the second and third terms in \eqref{eq_sec3_formula_for_integrand_1} cancel with the last term in \eqref{eq_sec3_formula_for_integrand_2}.
    Hence, \eqref{eq_int_int_dtW_G} can be rewritten as
    \begin{align*}
        & \int_0^\infty
            W_\varepsilon(z, t - \delta)
            G_\delta
            (
                \varepsilon
                + y
                - z
                + s(t)
                - s(t - \delta)
            )
        dz\\
        & - \int_0^\infty
            W_\varepsilon(z, {\delta})
            G_{t - {\delta}}
            (
                \varepsilon
                + y
                - z
                + s(t)
                - s(\delta)
            )
        dz\\
        & - \int_{\delta}^{t - \delta}
            f_\varepsilon(s(\tau), \tau)
            G_{t - \tau}
            (
                \varepsilon
                + y
                + s(t)
                - s(\tau)
            )
        d\tau\\
        & - \int_{\delta}^{t - \delta}
            \int_0^\infty
                \partial_z f_\varepsilon(s(\tau) + z, \tau)
                G_{t - \tau}
                (
                    \varepsilon
                    + y
                    - z
                    + s(t)
                    - s(\tau)
                )
            dz
        d\tau\\
        & +
        \int_{\delta}^{t - \delta}
            W_\varepsilon(0, \tau)
            \partial_y G_{t - \tau}
            (
                \varepsilon
                + y
                + s(t)
                - s(\tau)
            )
        d\tau.
    \end{align*}
    Therefore, the identity \eqref{eq_integration_wrt_ytau_boundary_ode} becomes
    \begin{align*}
            \int_\delta^{t - \delta}
                J_2
            d\tau
        & = - \int_{\delta}^{t - \delta}
            (
                \alpha
                - W_\varepsilon(0, \tau)
            )
            \dot{s}(\tau)
            G_{t - \tau}
            (
                \varepsilon
                + y
                + s(t)
                - s(\tau)
            )
        d\tau\\
        & - \int_{\delta}^{t - \delta}
                W_\varepsilon(0, \tau)
                (
                    \dot{S}_\varepsilon(\tau)
                    - \dot{s}(\tau)
                )
                G_{t - \tau}
                (
                    \varepsilon
                    + y
                    + s(t)
                    - s(\tau)
                )
            d\tau\\
        %& + W_\varepsilon(z, t)\\
        & + \int_0^\infty
            W_\varepsilon(z, t - \delta)
            G_\delta
            (
                \varepsilon
                + y
                - z
                + s(t)
                - s(t - \delta)
            )
        dz\\
        & - \int_0^\infty
            W_\varepsilon(z, {\delta})
            G_{t - \delta}
            (
                y
                - z
                + s(t)
                - s(\delta)
            )
        dz\\
        & + \int_{\delta}^{t - \delta}
            W_\varepsilon(0, \tau)
            \partial_y G_{t - \tau}
            (
                \varepsilon
                + y
                + s(t)
                - s(\tau)
            )
        d\tau\\
        & - \int_{\delta}^{t - \delta}
            f_\varepsilon(s(\tau), \tau)
            G_{t - \tau}
            (
                \varepsilon
                + y
                + s(t)
                - s(\tau)
            )
        d\tau\\
        & - \int_{\delta}^{t - \delta}
            \int_0^\infty
                \partial_z f_\varepsilon
                (
                    s(\tau) + z,
                    \tau
                )
                G_{t - \tau}
                (
                    \varepsilon
                    + y
                    - z
                    + s(t)
                    - s(\tau)
                )
            dz
        d\tau.
    \end{align*}
    First letting $\varepsilon \rightarrow 0$ and then $\delta \rightarrow 0$, we obtain from \eqref{eq_def_J1_J2}
    \begin{align*}
        0
        & = - \int_0^{t}
            (
                \alpha
                - w(0, \tau)
            ) \dot{s}(\tau) G_{t - \tau} (y + s(t) - s(\tau))
        d\tau\\
        & + w(y, t)\\
        & - \int_0^\infty
            w(z, 0) G_t(y - z + s(t))
        dz\\
        & + \int_0^{t}
            w(0, \tau) \partial_y G_{t - \tau} (y + s(t) - s(\tau))
        d\tau\\
        & - \int_0^{t}
            f(s(\tau), \tau) G_{t - \tau} (y + s(t) - s(\tau))
        d\tau\\
        & - \int_0^{t}
            \int_0^\infty
                \partial_z f(s(\tau); z, \tau)
                G_{t - \tau} (y - z + s(t) - s(\tau))
            dz
        d\tau.
    \end{align*}
    Note that $W_\varepsilon$ and $S_\varepsilon$ tend to $w$ and $s$, respectively, by uniqueness of the solution.
    Comparing the above identity with the solution formula of $w$, we obtain
    \begin{align*}
        0
        & = - \int_0^{t}
            (
                \alpha
                - w(0, \tau)
            ) \dot{s}(\tau) G_{t - \tau} (y + s(t) - s(\tau))
        d\tau\\
        & + \int_0^{t}
            (
                w(0, \tau)
                - \alpha
            ) \partial_y G_{t - \tau} (y + s(t) - s(\tau))
        d\tau.
    \end{align*}
    %for all $t>0$ and $y \in \Omega$.
    Then, taking the spatial trace, we find that
    \begin{align*}% \label{eq_formula_to_deduce_w_at_0_is_alpha}
        \int_0^{t}
            (
                w(0, \tau)
                - \alpha
            )
            \left[
                \dot{s}(\tau)
                - \frac{s(t) - s(\tau)}{2(t - \tau)}
            \right] G_{t - \tau}(s(t) - s(\tau))
        d\tau
        = 0
    \end{align*}
    for all $t > 0$.
    The part of the integrand following $w(0, \tau) - \alpha$ is positive and bounded from below, at least for short times, if $s$ is increasing.
    A similar argument applies if $s$ is decreasing.
    Therefore, by the monotonicity of $s$ and since $t$ is arbitrary, we conclude that
    \begin{align*}
        w(0, t)
        = \alpha
    \end{align*}
    for some sufficiently small $t_0 \in [0, T_\ast]$ and for a.a.\ $t \in [0, t_0]$.
    Repeating this argument with $t < 2t_0$, we obtain
    \begin{align*}% \label{eq_formula_to_deduce_w_at_0_is_alpha}
        \int_{t_0}^{t}
            (
                w(0, \tau)
                - \alpha
            )
            \left[
                \dot{s}(\tau)
                - \frac{s(t) - s(\tau)}{2(t - \tau)}
            \right] G_{t - \tau}(s(t) - s(\tau))
        d\tau
        = 0.
    \end{align*}
    Iteratively, we conclude that
    \begin{align*}
        w(0, t)
        = \alpha, \quad a.a. \; t \in [0, T_\ast].
    \end{align*}
\end{proof}

\section{Moving contact point condition} \label{sec_moving_boundary_condition}
This section is devoted to the monotonicity of $s$ obtained \red{in} Theorem \ref{thm_main1}.

\begin{lemma}
    Under the assumptions of Theorem \ref{thm_main1}, the solution $s$ is decreasing.
\end{lemma}

\begin{proof}
By Lemma \ref{lem_local_existence_of_w}\red{, we} can assume that $s$ is small in $X^{\frac12}(T)$.
Applying $\partial_y$ to \eqref{eq_def_of_w} and then the trace operator $\gamma_0$, together with \eqref{eq_formula_w_1_w_3}, we obtain
\begin{align} \label{eq_gamma_0_partial_y_w}
    \begin{split}
        & \gamma_0 \partial_y w(t)\\
        %& = - \int_0^\infty
        %    \partial_y G_t (z - s(t)) w_0(z)
        %dz\\
        %& + \int_0^t
        %    \partial_y G_{t - \tau} (s(t) - s(\tau)) f(s(\tau), \tau)
        %d\tau\\
        %& - \alpha G_{t} (s(t))
        %+ \alpha \int_0^t
        %    \dot{s}(\tau) \partial_y G_{t - \tau} (s(t) - s(\tau))
        %d\tau\\
        %& - \int_0^{t}
        %    \int_0^\infty
        %        \partial_y G_{t - \tau} (z - s(t) + s(\tau))
        %        \partial_z f(z + s(\tau), \tau)
        %    dz
        %d\tau\\
        & =
        %\frac{1}{2\sqrt{t}}
        %\int_{\frac{s(t)}{\sqrt{t}}}^\infty
        %    z
        %    e^{
        %        - z^2
        %    }
        %    \left(
        %        w_0
        %        (
        %            \sqrt{t} z
        %            - s(t)
        %        )
        %        - \alpha
        %    \right)
        %dz\\
        \int_0^\infty
            [\partial_y G_t](- z + s(t))
            \left(
                w_0(z)
                - \alpha
            \right)
        dz\\
        & + \int_0^t
            \gamma_0 \partial_y w(\tau) \frac{s(t) - s(\tau)}{2(t - \tau)} G_{t - \tau}(s(t) - s(\tau))
        d\tau\\
        & + \int_0^{t}
            \int_0^\infty
                \partial_y G_t(- z + s(t) - s(\tau))
                \partial_z f(z + s(\tau), \tau)
            dz
        d\tau.
    \end{split}
\end{align}
We define the integral operator $K$ by
\begin{align*}
    K \varphi
    = \int_0^t
        \frac{s(t) - s(\tau)}{2(t - \tau)} G_{t - \tau}(s(t) - s(\tau)) \varphi (\tau)
    d\tau.
\end{align*}
We estimate the singular behavior of the kernel as follows:
%\begin{align*}
%    & \left \vert
%        \int_0^t
%            \frac{
%                s(t) - s(\tau)
%            }{
%                4 \sqrt{\pi} (t - \tau)^{\frac{3}{2}} 
%            } e^{
%                - \frac{
%                    \vert s(t) - s(\tau) \vert^2
%                }{
%                    4 (t - \tau)
%                }
%            } \tau^{- \frac{1}{2}}
%        d\tau
%    \right \vert\\
%    & \leq \frac{1}{\sqrt{8 e \pi}} \left \vert
%        \int_0^{\eta t}
%            (t - \tau)^{- 1} \tau^{- \frac{1}{2}}
%        d\tau
%    \right \vert\\
%    & + \frac{1}{2 \sqrt{\pi}}
%    \Vert
%        s
%    \Vert_{X^\frac{1}{2}(T)}
%    \left \vert
%        \int_{\eta t}^t
%            (t - \tau)^{-\frac{1}{2}} \tau^{-1}
%        d\tau
%    \right \vert\\
%    & \leq \frac{
%        \eta^{\frac{1}{2}}
%    }{
%        \sqrt{2 e \pi}(1 - \eta)
%    } t^{- \frac{1}{2}}
%    + \frac{
%        \Vert
%            s
%        \Vert_{X^\frac{1}{2}(T)}
%        (1 - \eta)^\frac{1}{2}
%    }{
%        \sqrt{\pi} \eta
%    }
%    t^{- \frac{1}{2}}\\
%    & \leq \frac{
%        \Vert
%            s
%        \Vert_{X^\frac{1}{2}(T)}^\frac{1}{3}
%    }{
%        (2 e)^\frac{1}{6}
%    }
%    \left(
%        \frac{1}{(2e)^\frac{1}{3}}
%        + \Vert
%            s
%        \Vert_{X^\frac{1}{2}(T)}^\frac{2}{3}
%    \right)
%    t^{- \frac{1}{2}}
%    =: C_\ast t^{- \frac{1}{2}}.
%\end{align*}
\begin{align*}
    \left \vert
        \int_0^t
            \frac{
                s(t) - s(\tau)
            }{
                4 \sqrt{\pi} (t - \tau)^{\frac{3}{2}} 
            } e^{
                - \frac{
                    \vert s(t) - s(\tau) \vert^2
                }{
                    4 (t - \tau)
                }
            } \tau^{- \frac{1}{2}}
        d\tau
    \right \vert
    \leq C
    \Vert
        s
    \Vert_{X^\frac{1}{2}(T)}
    t^{- \frac{1}{2}}.
\end{align*}
By the construction of $s$, $C \Vert s \Vert_{X^\frac{1}{2}(T)}$ can be less than $\frac{1}{4}$.
We denote by $F$ the sum of the first and third terms \red{on} the right-hand side of \eqref{eq_gamma_0_partial_y_w}.
With this notation, $\gamma_0\partial_y w$ satisfies the Volterra integral equation of the second kind
\begin{align*}
    \gamma_0 \partial_y w(t)
    - K \gamma_0 \partial_y w(t)
    = F.
\end{align*}
%Since
%\begin{align*}
%    & \left \vert
%        \int_0^t
%            \frac{s(t) - s(\tau)}{2(t - \tau)} G_{t - \tau}(s(t) - s(\tau)) \tau^{- \frac{1}{2}}
%        d\tau
%    \right \vert\\
%    & \leq \frac{1}{\sqrt{8 e \pi}}
%    B \left(
%        \frac{1}{2}, \frac{1}{2}
%    \right)
%    = \sqrt{\frac{\pi}{8 e}}
%    =: C_\ast
%    < 1.
%\end{align*}
The Neumann series argument implies the existence of the inverse operator
\begin{align*}
    \gamma_0 \partial_y w
    = (I - K)^{-1} F
    = \sum_{k = 0}^{\infty} K^k F
\end{align*}
satisfying
\begin{align*}
    \sup_{0 < \tau < t} \tau^{\frac{1}{2}}
    \vert
        \gamma_0 \partial_y w(\tau)
    \vert
    \leq \frac{1}{
        1
        - C
        \Vert
            s
        \Vert_{X^\frac{1}{2}(T)}
    }
    \sup_{0 < \tau < t} \tau^{\frac{1}{2}}
    \vert
        F(\tau)
    \vert.
\end{align*}
Propositions \ref{prop_estimates_for_w1_w3} and \ref{prop_estimates_for_w4} imply that $t^\frac{1}{2} F(t)$ can be small for small $T>0$.
On the other hand, by assumption (C), $t^\frac{1}{2} f(s(t), t)$ is bounded from below.
This observation suggests that the principal contribution comes from $f(s(t),t)$, and then
\begin{align*}
    \alpha \dot{s}(t)
    = -f(s(t), t)
    - \gamma_0 \partial_y w(t)
    < 0
\end{align*}
for small $T_\ast < T$ and $t \in (0, T_\ast]$, \ie, $s$ is decreasing.
\end{proof}    

Summing up these \red{arguments}, we prove Theorem \ref{thm_main1}.
\begin{proof}[Proof of \red{Theorems} \ref{thm_main1} and \ref{thm_main2}.]
    The mild solution $(w, s)$ obtained in \red{Lemmas} \ref{lem_local_existence_of_w} and \ref{lem_local_existence_of_w_when_w0_f_are_smooth} satisfies $\dot{s}(t) < 0$ for $t \in (0, T)$ provided that $T$ is sufficiently small.
    Lemma \ref{lem_Neumann_boundary_condition} implies that such a solution $w$ satisfies the inhomogeneous Dirichlet boundary condition $\gamma_0 w = \alpha$.
    The function $v \in BUC([0, T); BUC(\Omega))$, defined by \eqref{eq_v_frome_w}, satisfies \eqref{eq_v} in the classical sense by construction of $w$, as stated in Remark \ref{rmk_boundary_condition}.
    Finally, by the maximum principle, we conclude that $v$ is positive on $\Omega$, and hence the front is located only at the origin.
\end{proof}

\section*{Acknowledgments}
The first author was partially supported by JSPS Grant-in-Aid for Young Scientists (No. JP26K17013) and by JST PRESTO, Japan, Grant Number JPMJPR26KB.
The work of the third author was partially supported by JSPS  through the grants KAKENHI:
JP24K00531, JP24H00183 and Daikin Industries, Ltd. through collaborative grants.

\end{document}